\documentclass[a4paper]{amsart}

\usepackage[english]{babel}
\usepackage[utf8x]{inputenc}
\usepackage[T1]{fontenc}
\usepackage{subdepth}
\usepackage{cite}

\usepackage{amsthm,amsmath,amssymb}
\usepackage{tikz-cd}
\usepackage{amsfonts}
\usepackage{graphicx}
\usepackage[colorlinks=true, allcolors=blue]{hyperref}
\usepackage{mathdots}
\usepackage{csquotes}
\usepackage{xfrac}

\theoremstyle{plain}
\newtheorem{thm}{Theorem}[section]
\newtheorem{lem}[thm]{Lemma}
\newtheorem{prop}[thm]{Proposition}
\newtheorem{cor}[thm]{Corollary}
\newtheorem*{thm*}{Theorem}
\numberwithin{equation}{section}
\newtheorem{thmx}{Theorem}

\theoremstyle{definition}
\newtheorem{example}[thm]{Example}
\newtheorem{defn}[thm]{Definition}

\theoremstyle{remark}

\newtheorem{rmk}[thm]{Remark}

\newcommand*{\suchthat}{\;\ifnum\currentgrouptype=16 \middle\fi|\;}
\newcommand{\Z}{\mathbb{Z}}

\DeclareMathOperator{\KK}{KK}

\newcommand{\acts}{\mathbin{\curvearrowright}}
\newcommand{\rightacts}{\mathbin{\curvearrowleft}}

\newcommand{\E}{\mathcal{E}}
\newcommand{\F}{\mathcal{F}}

\newcommand{\A}{\mathcal{A}}
\newcommand{\B}{\mathcal{B}}
\newcommand{\C}{\mathcal{C}}
\DeclareMathOperator{\Ab}{\mathbf{Ab}}

\newcommand{\bb}[1]{\mathbb{#1}}

\DeclareMathOperator{\Hom}{Hom}

\newcommand{\cal}[1]{\mathcal{#1}}
\newcommand{\bcdot}{\boldsymbol{\cdot}}
\renewcommand{\restriction}{\mathord{\upharpoonright}}

\DeclareMathOperator{\dom}{dom}
\DeclareMathOperator{\ran}{ran}

\DeclareMathOperator{\Ind}{Ind}

\DeclareMathOperator{\supp}{supp}

\newcommand{\hash}{\mathrel{\#}}

\newcommand{\defeq}{:=}
\newcommand{\orho}{\overline{\rho}}
\newcommand{\Cast}{\mathrm{C}^\ast}
\newcommand{\cs}{\Cast}

\newcommand\extrafootertext[1]{%
    \bgroup
    \renewcommand\thefootnote{\fnsymbol{footnote}}%
    \renewcommand\thempfootnote{\fnsymbol{mpfootnote}}%
    \footnotetext[0]{#1}%
    \egroup
}

\newenvironment{nalign}{
    \begin{equation}
    \begin{aligned}
}{
    \end{aligned}
    \end{equation}
    \ignorespacesafterend
}
\usepackage{environ}
\NewEnviron{killcontents}{}

\title[Functors between Kasparov categories from correspondences]{Functors between Kasparov categories from étale groupoid correspondences}
\author{Alistair Miller}
\address{Alistair Miller, Institut for Matematik og Datalogi, Syddansk Universitet, Campusvej 55, Odense M - DK-5230, Denmark, \href{mailto:mil@sdu.dk}{mil@sdu.dk}}

\begin{document}

\maketitle
\begin{abstract}
For an \'etale correspondence $\Omega \colon G \to H$ of \'etale groupoids, we construct an induction functor $\mathrm{Ind}_\Omega \colon \mathrm{KK}^H \to \mathrm{KK}^G$ between equivariant Kasparov categories. We introduce the crossed product of an $H$-equivariant correspondence by $\Omega$, and use this to build a natural transformation $\alpha_\Omega \colon K_*( G \ltimes \mathrm{Ind}_\Omega -) \Rightarrow K_*(H \ltimes -)$. When $\Omega$ is proper these constructions naturally sit above an induced map in K-theory $K_*(C^*(G)) \to K_*(C^*(H))$.
\end{abstract}

\extrafootertext{This work represents part of the author's PhD thesis, which was supported by the Engineering and Physical Sciences Research Council (EPSRC) through a doctoral studentship. The author has also been supported by the European Research Council (ERC) under the European Union's Horizon 2020 research and innovation programme (grant agreement No. 817597) and by the Independent Research Fund Denmark through the Grant 1054-00094B.}

\section{Introduction}

\'Etale groupoids provide models for a wide range of $\Cast$-algebras, encoding many kinds of topological dynamical system which occur naturally in examples. This range is demonstrated by Li's result that every classifiable $\Cast$-algebra has a (twisted) \'etale groupoid model \cite{Li20}. K-theory is a component of the Elliott invariant used to classify $\Cast$-algebras \cite{TWW17, EGLN15, Kirchberg95, Phillips00} and more generally extracts important information from a $\Cast$-algebra. The task of computing and controlling the K-theory of \'etale groupoid $\Cast$-algebras is the focus of a considerable programme of work, such as that surrounding the Baum--Connes conjecture \cite{Tu99b, BCH94, LeGall99, BoePro22}.

Controlling the K-theory as we move between different groupoids is a key technique for K-theory computations. For example, Morita equivalent groupoids have Morita equivalent $\Cast$-algebras \cite{MRW87} and thus isomorphic K-theory. This control is also a vital part of \cite{Li20}, where twisted groupoid $\Cast$-algebras with particular K-theory are constructed through an inductive system which behaves well with respect to K-theory. The connecting maps are realised at the groupoid level by a combination of open inclusions and certain proper homomorphisms in the opposite direction. Both these connecting maps and Morita equivalences are realised by groupoid correspondences. These are certain bispaces of groupoids and were introduced by Holkar (see \cite{Holkar14, Holkar17a, Holkar17b}) to serve as groupoid models for $\Cast$-correspondences, a kind of generalised morphism of $\Cast$-algebras used by Kasparov to build his bivariant K-theory. In particular, a proper correspondence induces a map in K-theory. \'Etale correspondences are covered in \cite{Albandik15, AKM22}, which place them as the morphisms in a (bi)category of \'etale groupoids and develop the functor from this (bi)category to the (bi)category of $\Cast$-algebras with $\Cast$-correspondences as morphisms.

The equivariant KK-theory of an \'etale groupoid $G$ \cite{LeGall99} is a framework through which we may better understand the K-theory groups $K_*(C^*(G))$ or $K_*(C^*_r(G))$. Working with the entire Kasparov category $\KK^G$ of $G$-$\Cast$-algebras and the K-theory functor $K_*(G \ltimes_r -) \colon \KK^G \to \Ab_*$ we may study $K_*(C^*_r(G)) = K_*(G \ltimes_r C_0(G^0))$ by relating $C_0(G^0)$ to other $G$-$\Cast$-algebras. This is seen in the categorical approach\cite{MeyNes06, BoePro22} to the Baum--Connes conjecture for $G$, which asserts that a certain approximation to $K_*(C^*_r(G))$ in terms of proper $G$-$\Cast$-algebras in $\KK^G$ is exact. The Baum--Connes conjecture has been verified for a wide range of \'etale groupoids, notably for amenable groupoids and more generally those with the Haagerup property \cite{Tu99b}. The ABC spectral sequence \cite{MeyNes10, Meyer08} was introduced to study Baum--Connes through the categorical approach. In the setting of a group $\Gamma$ satisfying the Baum--Connes conjecture this spectral sequence converges to the K-theory $K_*(C^*_r(\Gamma))$. For an \'etale groupoid $G$ with totally disconnected unit space and torsion-free isotropy groups satisfying Baum--Connes, an instance of the ABC spectral sequence connects the K-theory $K_*(C^*_r(G))$ to the groupoid homology $H_*(G)$ \cite{ProYam22}. This in turn has been used in low-dimensional situations \cite{BDGW23} to verify the HK conjecture \cite{Matui12} which computes $K_*(C^*_r(G))$ as the direct sum $\bigoplus_{n \in \bb Z} H_{* + 2n}(G)$.

Le Gall shows that Morita equivalent groupoids have equivalent Kasparov categories, and more generally constructs a functor between Kasparov categories from a generalised morphism of groupoids \cite[Section 7]{LeGall99}. B\"onicke uses this to construct an induction functor $\Ind^G_H \colon \KK^H \to \KK^G$ associated to a closed subgroupoid $H \subseteq G$ \cite[Proposition 4.7]{Boenicke20} in analogy with the group setting \cite[Section 5]{Kasparov95}, with the additional use of a forgetful functor $\KK^{G \ltimes X} \to \KK^G$. Groupoid correspondences are similar to the generalised morphisms Le Gall considers, more general in some ways and more specialised in others. There is a natural groupoid correspondence associated to a closed subgroupoid $H \subseteq G$. Our first main result takes a general \'etale correspondence $\Omega \colon G \to H$ and constructs a functor $\Ind_\Omega \colon \KK^H \to \KK^G$ based on B\"onicke's subgroupoid induction functor.
\begin{thmx}[Theorem \ref{KKinductionfunctor}]
A second countable Hausdorff \'etale correspondence $\Omega \colon G \to H$ induces an additive functor 
\[\Ind_\Omega \colon \KK^H \to \KK^G\]
of equivariant Kasparov categories which we call the induction functor.
\end{thmx}
The motivation for this stems from an effort to combine ideas surrounding the map in K-theory induced by a proper correspondence and those surrounding the Baum--Connes conjecture. In \cite{CEL13, CEL15, Li21a} there are isomorphisms in K-theory which may be expressed in the form
\begin{equation}\label{K theory formula}
 K_*(C^*_r(S)) \cong \bigoplus_{[e] \in E(S)^\times/S} K_*(C^*_r(S_e)) 
\end{equation}
for an inverse semigroup $S$ with stabiliser subgroups $S_e$ at non-zero idempotents $e \in E(S)^\times$. The method of proof requires that the dynamics of $S$ can be implemented by a group $\Gamma$. An equivariant Kasparov cycle in $\KK^\Gamma$ is then constructed, which provides the isomorphism (\ref{K theory formula}) by an application of the Going-Down principle\cite{CEOO04, Boenicke20}, which is a form of the Baum--Connes conjecture. For a general inverse semigroup $S$ there may be no such group $\Gamma$ but we might hope to replicate this procedure with the universal groupoid $S \ltimes \hat E$ of $S$. However, it does not seem possible to build the desired map (\ref{K theory formula}) in K-theory from a Kasparov cycle which is equivariant with respect to a single \'etale groupoid. Instead, it may be viewed as the map $K_*(C^*_r(G)) \to K_*(C^*_r(H))$ induced by an \'etale correspondence $\Omega \colon G \to H$. In order to apply powerful machinery connected to the Baum--Connes conjecture, we work to see how the map $K_*(C^*_r(G)) \to K_*(C^*_r(H))$ sits below behaviour induced by $\Omega$ at the higher level of the Kasparov categories $\KK^G$ and $\KK^H$.

The work in this document is part of a larger programme (see \cite{Miller24a, Miller24b}) which considers certain constructions associated to an \'etale groupoid $G$ and develops functoriality with respect to \'etale correspondences. Some of these constructions involve $\KK^G$, such as the spectral sequence $E(G)$ in \cite{ProYam22} linking $H_*(G)$ to $K_*(C^*_r(G))$ associated to an ample groupoid $G$ with torsion-free isotropy satisfying the Baum--Connes conjecture. In this programme we construct a morphism $E(G) \to E(H)$ of spectral sequences from a proper correspondence $\Omega \colon G \to H$ of such groupoids. This connects an induced map $H_*(G) \to H_*(H)$ to an induced map $K_*(C^*_r(G)) \to K_*(C^*_r(H))$, strengthening the link between groupoid homology and K-theory. Ultimately, K-theoretic isomorphisms follow from isomorphisms in groupoid homology. As an application, we deduce the K-theory isomorphism (\ref{K theory formula}) for each countable inverse semigroup $S$ whose universal groupoid $S \ltimes \hat E$ is Hausdorff, torsion-free and satisfies the Baum--Connes conjecture \cite{Miller24a}.

After the functoriality of the equivariant Kasparov category $\KK^G$ associated to an \'etale groupoid $G$ is established, the next step towards understanding the behaviour of $K_*(C^*_r(G))$ is the behaviour of the K-theory functor $K_*(G \ltimes_r -) \colon \KK^G \to \Ab_*$. From an \'etale correspondence $\Omega \colon G \to H$ we connect the K-theory functors $K_*(G \ltimes -) \colon \KK^G \to \Ab_*$ and $K_*(H \ltimes -) \colon \KK^H \to \Ab_*$. We use the full crossed product as it has greater functoriality properties than the reduced crossed product\footnote{A mild amenability condition may be imposed on the correspondence in order to work at the reduced level, and through the Baum--Connes conjecture we can understand the reduced crossed product by working with the full one.}. In order to obtain a natural transformation we develop a version of the $\Cast$-correspondence construction $C^*(\Omega) \colon C^*(G) \to C^*(H)$ \cite{AKM22} which allows for coefficients in $\Cast$-algebras\footnote{It is possible to treat a more general kind of $\Cast$-coefficients, see \cite{Miller22}, but we do not have any applications for this in mind.}.
\begin{thmx}[Proposition \ref{crossed product correspondence}]
A Hausdorff \'etale correspondence $\Omega \colon G \to H$ and an $H$-equivariant $\cs$-correspondence $E \colon A \to B$ induce a $\cs$-correspondence
\[\Omega \ltimes E \colon G \ltimes \Ind_\Omega A \to H \ltimes B\] 
which we call the crossed product of $E$ by $\Omega$.
\end{thmx}
The crossed product correspondence $\Omega \ltimes E$ is proper when $E \colon A \to B$ is proper. In particular if we take an identity correspondence $A \colon A \to A$ then $\Omega \ltimes A$ induces a map in K-theory from $ K_*(G \ltimes \Ind_\Omega A)$ to $K_*(H \ltimes A)$.
\begin{thmx}[Theorem \ref{induction natural transformation KK}]
Let $\Omega \colon G \to H$ be a second countable Hausdorff \'etale correspondence. Then the assignment of the map
\[ K_*(\Omega \ltimes A) \colon K_*(G \ltimes \Ind_\Omega A) \to K_*(H \ltimes A) \]
to an $H$-$\Cast$-algebra $A \in \KK^H$ is natural with respect to morphisms in $\KK^H$, so determines a natural transformation 
\[ \alpha_\Omega \colon K_*(G \ltimes \Ind_\Omega -) \Rightarrow K_*(H \ltimes -) \colon \KK^H \rightrightarrows \Ab_*. \]
\end{thmx}
Specialising to the $H$-$\Cast$-algebra $C_0(H^0)$, the induced algebra $\Ind_\Omega C_0(H^0)$ is given by $C_0(\Omega/H)$ and we obtain a map $K_*(C^*(G \ltimes \Omega/H)) \to K_*(C^*(H))$. We cannot in general get a map from $K_*(C^*(G))$ unless we assume that $\Omega \colon G \to H$ is proper. In this case, there is a $G$-equivariant proper map $\orho \colon \Omega/H \to G^0$ which induces an element $f_\Omega \in \KK^G(C_0(G^0),C_0(\Omega/H))$. This enables us to recover the map induced in K-theory by $\Omega$ as the composition
\[ \alpha_\Omega(C_0(H^0)) \circ K_*(G \ltimes f_\Omega) \colon K_*(C^*(G)) \to K_*(C^*(H)).\]
We also verify that all our constructions are compatible with the composition of \'etale correspondences. As isomorphisms in this category are given by Morita equivalences, this in particular recovers the Morita invariance of the equivariant Kasparov category $\KK^G$ and the K-theory functor $K_*(G \ltimes -) \colon \KK^G \to \Ab_*$ in a very concrete fashion.

In Section \ref{preliminaries}, we review material which we use to work concretely in groupoid equivariant KK-theory. This includes the introduction of a strict topology for bundles of adjointable operators for Hilbert modules over a $C_0(X)$-algebra. In Section \ref{correspondence section} we give examples of \'etale correspondences, including new examples from inverse semigroups. We develop the induced Banach space construction associated to an \'etale correspondence in Section \ref{induced banach space section} and the crossed product correspondence construction in Section \ref{crossed product section}. In Section \ref{KK section} these are used in turn to build the KK-theoretic induction functor and the natural transformation associated to an \'etale correspondence. Finally, compatibility with composition of correspondences is verified in Section \ref{compatibility with composition section}.

The author would like to thank Christian B\"onicke, Xin Li and Ralf Meyer for helpful comments and discussions.

\section{Preliminaries}\label{preliminaries}

\subsection{\'Etale groupoids, inverse semigroups and their actions}

A \textit{topological groupoid} $G$ is a topological space $G$ with a distinguished subspace $G^0$ called the \textit{unit space}, with continuous maps $r,s \colon G \rightrightarrows G^0$ that assign a \textit{range} and a \textit{source} to each groupoid element. Elements $g$ and $h$ with matching range and source $r(h) = s(g)$ can be \textit{multiplied} or \textit{composed} to form an element $gh$ with range $r(g)$ and source $s(h)$. Multiplication is a continuous associative map $(g,h) \mapsto gh \colon G^2 \to G$, where $G^2 \defeq \{(g,h) \in G \times G \suchthat s(g) = r(h) \}$ is the space of composable pairs. Elements of $G^0$ act as identities in that each element is its own range and source, and $r(g) g = g = g s(g)$ for each $g \in G$. There is a continuous \textit{inversion} map $g \mapsto g^{-1} \colon G \to G$ that swaps the range and source of each element such that $g^{-1} g = s(g)$ and $g g^{-1} = r(g)$ for each $g \in G$. 

For each $x \in G^0$ we write $G^x$ for the range fibre $r^{-1}(x)$ and $G_x$ for the source fibre $s^{-1}(x)$ and their intersection is the \textit{isotropy group} $G^x_x$. The topological groupoid $G$ is \textit{\'etale} if the range map $r \colon G \to G^0$ is a local homeomorphism and the unit space $G^0$ is open in $G$. In this setting, the source map is also a local homeomorphism, and a subset $U \subseteq G$ is a \textit{bisection} if $r$ and $s$ are injective on $U$. Unless otherwise stated, we assume that the unit space $G^0$ of an \'etale groupoid is locally compact and Hausdorff, in which case $G$ is locally compact and locally Hausdorff, but not necessarily Hausdorff.

An \textit{inverse semigroup} $S$ is a semigroup such that for each $s \in S$ there is a unique element $s^* \in S$ called the \textit{adjoint} of $s$ satisfying $s^*ss^* = s^*$ and $ss^*s=s$, and by convention we always include a distinguished $0$ element satisfying $0 s = 0$ and $s 0 = 0$ for each $s \in S$. The set $E(S)$ of idempotents forms a semilattice under multiplication. 

We think of the elements of an inverse semigroup as partial symmetries, and this is reflected in how they act. An inverse semigroup acts on a space $X$ by \textit{partial homeomorphisms}, which are homeomorphisms $\alpha \colon \dom \alpha \to \ran \alpha$ between open subsets of $X$. The composition $\alpha \circ \beta$ of two partial homeomorphisms is the partial homeomorphism $x \mapsto \alpha(\beta(x))$ whose domain consists of the $x$ in $\dom \beta$ with $\beta(x) \in \dom \alpha$. The set of partial homeomorphisms of $X$ forms an inverse semigroup $\text{PHom}(X)$, with the adjoint $\alpha^*$ of a partial homeomorphism $\alpha$ given by the inverse $\alpha^{-1} \colon \ran \alpha \to \dom \alpha$.

\begin{defn}[Inverse semigroup action]
A \textit{(left) action} $S \acts X$ of an inverse semigroup $S$ on a space $X$ is a homomorphism $\alpha \colon S \to \text{PHom}(X)$ to the inverse semigroup of partial homeomorphisms of $X$. The domain and range of $\alpha(s)$ for $s \in S$ may be written $\dom_X s$ and $\ran_X s$ if we want to emphasise the space. We usually write $s \bcdot x$ for $\alpha(s)(x)$.
\end{defn}

A continuous map $f \colon X \to Y$ between $S$-spaces is \textit{$S$-equivariant} if $x \in \dom_X s$ implies that $f(x) \in \dom_Y s$ and $f(s \bcdot x) = s \bcdot f(x)$. From an inverse semigroup action we may build an \'etale groupoid:

\begin{defn}[Transformation groupoid of an inverse semigroup action]
Let $S \acts X$ be an action of an inverse semigroup $S$ on a locally compact Hausdorff space $X$. Consider the quotient space
\begin{equation*}
S \ltimes X \defeq \{ (s,x) \in S \times X \suchthat x \in \dom s \}/\sim,
\end{equation*}
where $(s,x) \sim (t,x)$ if there is an idempotent $e \in E$ with $x \in \dom e$ and $se = te$. We take the quotient topology of the subspace topology of the product topology on $S \times X$, where $S$ is given the discrete topology, and we write $[s,x]$ for the equivalence class of $(s,x)$. We define range and source maps by $r([t,x]) = t \bcdot x$, $s([t,x]) = x$ for $[t,x] \in S \ltimes X$ and a composition by $[t,y] [s,x] = [ts,x]$ when $y = s \bcdot x$. This turns $S \ltimes X$ into a (possibly non-Hausdorff) \'etale groupoid called the \textit{transformation groupoid} of $S \acts X$ whose unit space is an open subset of $X$. Each \'etale groupoid may be obtained this way, for example through the action of its inverse semigroup of Hausdorff open bisections on its unit space.
\end{defn}

\begin{defn}[Groupoid action]
Let $G$ be an \'etale groupoid and let $X$ be a topological space. A \textit{(left) groupoid action} $G \acts X$ consists of:
\begin{itemize}
\item a continuous map $\tau \colon X \to G^0$ called the \textit{anchor map},
\item a continuous map $\alpha \colon G \times_{s,G^0,\tau} X \to X$ called the \textit{action map}. We usually write $g \bcdot x$ for $\alpha(g,x)$.
\end{itemize}
The pair $(\tau, \alpha)$ is an action if whenever $g, h \in G$ and $x \in X$ satisfy $s(g) = r(h)$ and $s(h) = \tau(x)$, we have $\tau(h \bcdot x) = r(h)$ and $gh \bcdot x = g \bcdot (h \bcdot x)$. Each groupoid element $g \in G$ induces a homeomorphism $\alpha_g \colon X_{s(g)} \to X_{r(g)}$ that sends $x$ to $g \bcdot x$, where we write $X_z = \tau^{-1}(z)$ for the fibre of $X$ at $z \in G^0$. We call $X$ a \textit{$G$-space}. 
\end{defn}

There is a canonical action of an \'etale groupoid $G$ on its unit space $G^0$ for which $g \bcdot s(g) = r(g)$ for each $g \in G$. This allows us to view $G$ through the lens of an action. Conversely, we may view an action $G \acts X$ through the lens of a groupoid via the action groupoid $G \ltimes X$:

\begin{defn}[Action groupoid]
Let $G$ be an \'etale groupoid and let $X$ be a locally compact Hausdorff $G$-space with anchor map $\tau \colon X \to G^0$. The fibre product $G \times_{s,G^0,\tau} X$ may be given the structure of an \'etale groupoid $G \ltimes X$ with unit space $X$ which we call the \textit{action groupoid}. For $(g,x) \in G \ltimes X$ the range is given by $g \bcdot x$ and the source by $x$, and multiplication is given by $(g,x) (h,y) = (gh,y)$ for $(g,x),(h,y) \in G \ltimes X$ with $x = h \bcdot y$.
\end{defn}

\begin{defn}[Free, proper and \'etale actions]
Let $G$ be an \'etale groupoid and let $X$ be a $G$-space with anchor map $\tau \colon X \to G^0$. We say that the $G$-space $X$ or the action $G \acts X$ is:
\begin{itemize}
\item \textit{free} if for $g \in G$ and $x \in X$,  $g \bcdot x = x$ implies that $g \in G^0$.
\item \textit{proper} if the map $ (g,x) \mapsto (g \bcdot x, x) \colon G \times_{s,G^0,\tau} X \to X \times X$ is proper.
\item \textit{\'etale} if $\tau \colon X \to G^0$ is a local homeomorphism.
\end{itemize}
We say that $G$ is \textit{proper} if the canonical action $G \acts G^0$ on its unit space is proper, and we say that $G$ is \textit{principal} if this action is free.
\end{defn}

We collect a number of key facts about proper actions.

\begin{prop}[Proposition 2.12 in \cite{Tu04}, Lemma 2.12 and Proposition 2.19 in \cite{AKM22}]\label{free proper action properties}
Let $G$ be an \'etale groupoid and let $X$ be a proper $G$-space. Then the orbit space $G \backslash X$ is Hausdorff. If the action is also free, then the quotient map $q \colon X \to G \backslash X$ is a local homeomorphism.
\end{prop}

\begin{prop}[Proposition 2.20 in \cite{Tu04}]
Let $G$ be an \'etale groupoid, let $X$ be a Hausdorff right $G$-space and let $Y$ be a proper left $G$-space. Then the diagonal action $X \times_{G^0} Y \rightacts G$ is proper.
\end{prop}

\begin{defn}[Pullback over a groupoid]
Let $G$ be an \'etale groupoid, let $X$ be a right $G$-space with anchor map $\tau$ and let $Y$ be a left $G$-space. Consider the diagonal action of $G$ on $X \times_{G^0} Y$ given by $(x,y) \bcdot g := (x \bcdot g, g^{-1} \bcdot y)$ for $ (x,y) \in X \times_{G^0} Y$ and $g \in G^{\tau(x)}$. The \textit{pullback} $X \times_G Y$ of $X$ and $Y$ over $G$ is the quotient space $(X \times_{G^0} Y)/G$ of the diagonal action.
\end{defn}

\subsection{$C_0(X)$-Banach spaces and Banach bundles}

\begin{defn}[$C_0(X)$-Banach space]
Let $X$ be a locally compact Hausdorff space. A $C_0(X)$-Banach space is a Banach space $A$ equipped with a bilinear map $(f,a) \mapsto f \bcdot a \colon C_0(X) \times A \to A$ called the \textit{structure map} such that
\begin{itemize}
\item for each $f \in C_0(X)$ and $a \in A$, we have $\lVert f \bcdot a \rVert \leq \lVert f \rVert \lVert a \rVert$, and
\item for each $f,g \in C_0(X)$ and $a \in A$, we have $f \bcdot (g \bcdot a) = fg \bcdot a$, and
\item elements of the form $f \bcdot a$ for $f \in C_0(X)$ and $a \in A$ have dense span in $A$.
\end{itemize}
In short, $C_0(X)$ acts non-degenerately on $A$. We usually omit the $\bcdot$ and write $fa$ instead of $f \bcdot a$. A continuous linear map $\varphi \colon A \to B$ of $C_0(X)$-Banach spaces is \textit{$C_0(X)$-linear} if $f \varphi(a) = \varphi(f  a)$ for each $f \in C_0(X)$ and $a \in A$. 
\end{defn}
A $\cs$-algebra $B$ is called a $C_0(X)$-algebra if it is equipped with a non-degenerate $*$-homomorphism $C_0(X) \to ZM(B)$ into the central multipliers of $B$. This gives $B$ the structure of a $C_0(X)$-Banach space, and moreover if $E$ and $F$ are Hilbert $B$-modules, each of $E$, $F$ and $\cal K(E,F)$ are $C_0(X)$-Banach spaces (see \cite[Section 4]{Boenicke20}), but the action of $C_0(X)$ on $\cal L(E,F)$ may fail to be non-degenerate. The \textit{fibre} $A_x$ at $x \in X$ of a $C_0(X)$-Banach space $A$ is the quotient of $A$ by the closed subspace $\overline{C_0(X \setminus \{x\}) A}$. Each $a \in A$ determines a section $x \mapsto a_x \colon X \to \bigsqcup_{x \in X} A_x$. As for $C_0(X)$-algebras, the function $x \mapsto \lVert a_x \rVert \colon X \to \bb R_{\geq 0}$ is upper semicontinuous and vanishes at infinity, and for each $f \in C_0(X)$ and $x \in X$ we have $(fa)_x = f(x)a_x$.
We say a $C_0(X)$-Banach space $A$ is \textit{$C_0(X)$-locally convex} if $\lVert a \rVert = \sup_{x \in X} \lVert a_x \rVert$ for each $a \in A$. Every $C_0(X)$-algebra $B$ is locally $C_0(X)$-convex \cite[Proposition C.10]{Williams07}. Each fibre $E_x$ of a Hilbert $B$-module $E$ is a Hilbert $B_x$-module with operations compatible with the fibre map $E \to E_x$ \cite[Section 4]{Boenicke20}, and local $C_0(X)$-convexity follows from that of $B$. Given another Hilbert $B$-module $F$, each $T \in \cal L(E, F)$ is automatically $C_0(X)$-linear so there is a unique $T_x \in \cal L(E_x,F_x)$ such that $(Te)_x = T_x e_x$ for each $e \in E$. Moreover, we have $\lVert T \rVert = \sup_{x \in X} \lVert T_x \rVert$ (see \cite[Lemma C.11]{Williams07}). If $T$ is compact then so is $T_x$, and this induces an isomorphism $\cal K(E,F)_x \cong \cal K(E_x,F_x)$ identifying the fibre of $T$ at $x$ with $T_x$\footnote{This appears as \cite[Corollaire 2.3.1]{LeGall94} and \cite[Proposition 4.2]{LeGall99}, and while the proof is difficult to track down, its main difficulty is the injectivity. This follows from the fact that $x \mapsto \lVert T_x \rVert$ is upper semicontinuous and vanishes at infinity for each $T \in \cal K(E,F)$. The upper semicontinuity follows in turn from finding $b \in M_n(B)$ for each finite rank $T$ such that $\lVert T_x \rVert = \lVert b_x \rVert$ for each $x \in X$, see \cite[Proposition 1.51]{Miller22} for details.}. To understand these continuous families of fibres we turn to Banach bundles:

\begin{defn}[Banach bundle]
A \textit{Banach bundle} over a locally compact Hausdorff space $X$ is a topological space $\cal A$ equipped with a continuous, open surjection $p \colon  \cal A \to X$ and complex Banach space structures on each fibre $A_x := p^{-1} (x)$ such that
\begin{itemize}
\item the map $a \mapsto \lVert a \rVert \colon \A \to \bb R_{\geq 0}$ is upper semicontinuous,
\item the map $(a,b) \mapsto a + b \colon \A \times_X \A \to \A$ is continuous,
\item for each $\lambda \in \bb C$, the map $a \mapsto \lambda a \colon \A \to \A$ is continuous, and
\item if $(a_i)_i$ is a net in $\cal A$ such that $p(a_i) \to x$ and $\lVert a_i \rVert \to 0$, then $a_i \to 0_x \in A_x$.
\end{itemize}
A \textit{map} $\varphi \colon \A \to \B$ of Banach bundles over $X$ is a function such that $\varphi(A_x) \subseteq B_x$ for each $x \in X$ and these induced maps $(\varphi_x \colon A_x \to B_x)_{x \in X}$ form a bounded collection of continuous linear maps. We are interested in the continuous maps of Banach bundles. A \textit{$\Cast$-bundle} is a Banach bundle $\A \to X$ with the structure of a $\Cast$-algebra on each fibre such that multiplication $\A \times_X \A \to \A$ and the adjoint map $\A \to \A$ are continuous.
\end{defn}

We may obtain Banach bundles from $C_0(X)$-Banach spaces through the following construction, see II.13.18 in \cite{DorFel88a} and Theorem C.25 in \cite{Williams07}. While \cite{DorFel88a} concerns only Banach bundles with a continuous norm function and \cite{Williams07} concerns only $\Cast$-bundles, their proofs apply equally well to the general setting of Banach bundles.

\begin{prop}\label{Banach bundle topology construction}
Let $X$ be a locally compact Hausdorff space, let $(A_x)_{ x \in X }$ be a collection of Banach spaces and let $\Gamma$ be a complex vector space of sections $X \to \cal A = \bigsqcup_{x \in X} A_x$ such that
\begin{itemize}
\item for each $a \in \Gamma$, $x \mapsto \lVert a(x) \rVert$ is upper semicontinuous,
\item for each $x \in X$, the set $\{a(x) \suchthat  a \in \Gamma \}$ is dense in $A_x$.
\end{itemize} 
Then there is a unique topology on $\cal A$ turning $\cal A \to X$ into a Banach bundle such that each section in $\Gamma$ is continuous. 
\end{prop}
Therefore for any $C_0(X)$-Banach space $A$, there is a unique Banach bundle structure on $\A = \bigsqcup_{x \in X} A_x \to X$ such that $x \mapsto a_x \colon X \to \A$ is a continuous section for each $a \in A$. We call this the Banach bundle associated to the $C_0(X)$-Banach space $A$, and by default we refer to it using the calligraphic version of the same letter.

Conversely, for any Banach bundle $\A \to X$, the space $\Gamma_0(X, \A)$ of continuous sections which vanish at infinity is a $C_0(X)$-Banach space with the pointwise action of $C_0(X)$. For each $x \in X$ the evaluation at $x$ map $\Gamma_0(X, \A) \to A_x$ identifies the fibre of $\Gamma_0(X,\A)$ at $x$ with $A_x$. The key step in the proof of this is surjectivity of the evaluation map, see Theorem C17 and Remark C18 in \cite{DorFel88a}. Proposition \ref{Banach bundle topology construction} recovers $\A$ from $\Gamma_0(X,\A)$ as the associated Banach bundle.

For any $C_0(X)$-Banach space $A$, the section $x \mapsto a_x \colon X \to \A$ associated to each element $a \in A$ determines a $C_0(X)$-linear contraction $A \to \Gamma_0(X, \A)$. This is isometric exactly when $A$ is locally $C_0(X)$-convex, as is each $C_0(X)$-Banach space we utilise for equivariant KK-theory. The map $A \to \Gamma_0(X, \A)$ always has dense image by the following density criterion, and thus is an isomorphism when $A$ is locally $C_0(X)$-convex:

\begin{prop}[Proposition C.24\footnote{The proof carries over word for word to general Banach bundles.} in \cite{Williams07}]\label{density in Banach spaces}
Let $\A \to X$ be a Banach bundle over a locally compact Hausdorff space $X$ and let $\Gamma \subseteq \Gamma_0(X, \A)$. Suppose that $\Gamma$ is closed under the action of $C_0(X)$ and for each $x \in X$, the set $\{ \gamma_x \suchthat \gamma \in \Gamma \}$ has dense span in $A_x$. Then $\Gamma$ has dense span in $\Gamma_0(X,\A)$.
\end{prop}

We thus frequently identify the elements $a \in A$ of a locally $C_0(X)$-convex $C_0(X)$-Banach space with their associated sections $x \mapsto a_x \in \Gamma_0(X,\A)$, and switch freely between studying the $C_0(X)$-Banach space $A$ and the Banach bundle $\A \to X$. We collect a number of results concerning the topology of Banach bundles.

\begin{lem}[Proposition C.20 in \cite{Williams07}]\label{topology of banach bundle}
Let $p \colon \A \to X$ be a Banach bundle over a locally compact Hausdorff space $X$ and let $( a_i )_i$ be a net in $\A$ such that $p(a_i) \to p(a)$ for some $a \in \A$. Suppose that for all $\epsilon > 0$ there is a net $(u_i)_i$ in $\A$ and $u \in \A$ with $p(u_i) = p(a_i)$ and $p(u) = p(a)$ such that
\begin{itemize}
\item $u_i \to u$ in $\A$,
\item $\lVert a - u \rVert < \epsilon$, and
\item $\lVert a_i - u_i \rVert < \epsilon$ for large $i$.
\end{itemize}
Then $a_i \to a$.
\end{lem}

Straightforward applications of this lemma yield conditions for the continuity of maps and bilinear maps of Banach bundles:
\begin{prop}[Continuity condition for maps out of bundles]\label{mapoutofbundle}
Let $\varphi \colon  \A \to \B$ be a map of Banach bundles over a locally compact Hausdorff space $X$ and let $\Gamma \subseteq \Gamma_0(X, \A)$ be a subset with dense span. If $\varphi \circ \gamma$ is continuous for each $\gamma \in \Gamma$ then $\varphi$ is continuous. Similarly, let $\psi \colon \A \times_X \B \to \C$ be a bounded collection of bilinear or sesquilinear maps for Banach bundles $\A$, $\B$ and $\C$, and suppose $\Gamma_A \subseteq \Gamma_0(X,\A)$ and $\Gamma_B \subseteq \Gamma_0(X,\B)$ have dense span. If $\psi \circ (\gamma_A, \gamma_B) \colon X \to \C$ is continuous for each $\gamma_A \in \Gamma_A$ and $\gamma_B \in \Gamma_B$ then $\psi$ is continuous.
\end{prop}

\begin{rmk}\label{bilinear continuity remark}
This continuity condition allows us to check that the fibrewise multiplication and adjoint on the bundle associated to a $C_0(X)$-algebra $B$ are continuous. Moreover, given Hilbert $B$-modules $E$ and $F$, the fibrewise module and inner product structures define continuous maps, and the map $\Theta_{-,-} \colon \F \times_X \E \to \cal K(\E,\F)$ is continuous.
\end{rmk}

The other continuity lemma we will utilise is the following:

\begin{lem}[Condition for isomorphism of Banach bundles, see Remark 3.8 in \cite{MuhWil08a}]\label{Banach bundle isomorphism condition}
Let $\A \to X$ and $\B \to X$ be Banach bundles over a locally compact Hausdorff space $X$ and let $\varphi \colon  \A \to \B$ be a continuous map of Banach bundles such that $\varphi_x \colon  A_x \to B_x$ is an isometric isomorphism for each $x \in X$. Then $\varphi$ is an isomorphism of Banach bundles.
\end{lem}

We might hope that given Hilbert modules $E$ and $F$ over a $C_0(X)$-algebra, there is a Banach bundle structure on the space $\cal L(\E,\F) = \bigsqcup_{x \in X} \cal L(E_x,F_x)$ such that the assignment mapping $T \in \cal L(E,F)$ to the section $x \mapsto T_x \colon X \to \cal L(\E,\F)$ is an isomorphism $\cal L(E,F) \cong \Gamma_b(X, \cal L(\E,\F))$. Unfortunately this is not always possible, as $x \mapsto \lVert T_x \rVert \colon X \to \bb R$ may fail to be upper semicontinuous \cite[Remark C.14]{Williams07}. Nevertheless, it is still useful to consider a topology on the total space $\cal L(\E,\F)$. As in \cite[Lemma C.11]{Williams07}, a bounded section $x \mapsto T_x \colon X \to \cal L(\E,\F)$ determines an adjointable operator $T \in \cal L(E,F)$ if and only if for each $e \in E$ and $f \in F$ the sections $x \mapsto T_x e_x$ and $x \mapsto T^*_x f_x$ are continuous (and therefore identified with elements of $F$ and $E$ respectively). This may be viewed naturally as a continuity condition on $x \mapsto T_x$ for a topology on $\cal L(\E,\F)$.

\begin{defn}[Strict topology for bundles of adjointable operators]
Let $E$ and $F$ be Hilbert $B$-modules for a $C_0(X)$-algebra $B$ and let $\cal L(\E,\F) = \bigsqcup_{x \in X} \cal L(E_x,F_x)$. The \textit{strict topology} on $\cal L(\E,\F)$ is the weakest topology such that the maps
\begin{align*}
\cal L(\E,\F) & \to \F & \cal L(\E,\F) & \to \E \\
(x,T) & \mapsto (x,T(e_x)) & (x,T) & \mapsto (x,T^*(f_x))
\end{align*}
are continuous for each $e \in E$ and $f \in F$. 
\end{defn} 
We may then identify an operator $T \in \cal L(E,F)$ with its bounded strictly continuous section $x \mapsto T_x \in \Gamma_b(X, \cal L(\E,\F))$. A bounded section $x \mapsto T_x$ is strictly continuous if and only if the maps $(T_x)_{x \in X} \colon \E \to \F$ and $(T^*_x)_{x \in X} \colon \F \to \E$ of Banach bundles are continuous.

\subsection{Actions of groupoids on Banach spaces}

Given a continuous map $f \colon Y \to X$ of locally compact Hausdorff spaces and a Banach bundle $\A \to X$, the pullback bundle $f^* \A \defeq Y \times_X \A \to Y$ is a Banach bundle whose fibre at $y \in Y$ may be identified with $A_{f(y)}$. We may then define the pullback $f^* A$ of a locally $C_0(X)$-convex $C_0(X)$-Banach space $A$ to be the section space $\Gamma_0(Y, f^* \A)$. The pullback $f^* B$ of a $C_0(X)$-algebra $B$ by $f \colon Y \to X$ is a $C_0(Y)$-algebra, and the pullback $f^* E$ of a Hilbert $B$-module $E$ is a Hilbert $f^* B$-module with operations given pointwise over $Y$. Given another Hilbert $B$-module $F$, we may identify $f^* \cal K(E,F)$ with $\cal K(f^* E, f^* F)$.
\begin{rmk}\label{pullback section remark}
A section $a \colon Y \to f^* \A$ can be identified with a function $a \colon Y \to \A$ such that $a(y) \in A_{f(y)}$ for each $y \in Y$, and we will frequently make this identification. Continuity of one implies continuity of the other.
\end{rmk}

\begin{rmk}\label{pullbacks allow us to describe restrictions}
The restriction $\A \restriction_Y \to Y$ of a Banach bundle $\A \to X$ to a subspace $Y \subseteq X$ may be described as the pullback of $\A \to X$ with respect to the inclusion $Y \hookrightarrow X$. For the section space we write $\Gamma_0(Y, \A)$ instead of $\Gamma_0(Y, \A \restriction_Y)$.
\end{rmk}

For a continuous map $f \colon Y \to X$ of locally compact Hausdorff spaces and a continuous section $\xi \colon Y \to \A$ to a Banach bundle $\A \to Y$, we say $\xi$ has \textit{proper support} with respect to $f$ if the restriction of $f$ to the (closed) support $\supp \xi$ is proper. The following version of \cite[Lemma 3.12]{Boenicke20} has a near identical proof.

\begin{lem}[Sum over fibres lemma]\label{propersupportlemma}
Let $f\colon X \to Y$ be a local homeomorphism of locally compact Hausdorff spaces and let $\A \to Y$ be a Banach bundle. Let $\xi \in \Gamma_b(X, f^* \A)$ be a bounded continuous section with proper support with respect to $f$. Then the section $f_* \xi \colon  Y \to \A$ given by
\[ f_* \xi (y) = \sum_{x \in f^{-1}(y)} \xi(x) \]
is well-defined and continuous. If $\xi$ is compactly supported then so is $f_* \xi$.
%
\end{lem}

\begin{defn}[Groupoid action on a Banach space]
Let $G$ be a Hausdorff \'etale groupoid with unit space $X$, and let $A$ be a locally $C_0(X)$-convex $C_0(X)$-Banach space. A \textit{(left) Banach space action} $\alpha \colon G \acts A$ of $G$ on $A$ is an isometric isomorphism $\alpha \colon s^* \A \to r^* \A$ of Banach bundles such that
\begin{itemize}
\item for each $x \in X$ and $a \in A_x$, we have $\alpha(x,a) = a$, and
\item for each $(g,h) \in G^2$ and $a \in A_{s(h)}$, we have $\alpha(g, \alpha(h,a)) = \alpha(gh,a)$.
\end{itemize}
In other words, this is a continuous action of $G$ on $\A$ by linear isometric isomorphisms. We usually suppress $\alpha$ and write $g \bcdot a$ for $\alpha(g,a)$. We call $A$ a \textit{$G$-Banach space}\footnote{Note that local $C_0(X)$-convexity is now implicit, as these are the only Banach spaces we deal with.} and $\A$ a \textit{Banach $G$-bundle}. 
\end{defn}
It is often convenient to check that an action map is continuous with the following lemma, see \cite[Lemma 3.9]{Boenicke20}:
\begin{lem}\label{continuity of action map}
Let $G$ be a Hausdorff \'etale groupoid, let $\A$ be a Banach $G$-bundle and suppose $\Gamma \subseteq \Gamma_0(G^0,\A)$ has dense span. Then a map $(g,a) \mapsto g \bcdot a \colon s^* \A \to r^* \A$ of Banach bundles is continuous if for each $\xi \in \Gamma$ the map $g \mapsto g \bcdot \xi_{s(g)} \colon G \to r^*\A$ is continuous.
\end{lem}

A Banach space action of a Hausdorff \'etale groupoid $G$ on a $C_0(G^0)$-algebra $B$ is a \textit{$\cs$-algebra action} if for each $g \in G$ the induced map $B_{s(g)} \to B_{r(g)}$ is a $*$-isomorphism, in which case $B$ is a \textit{$G$-$\cs$-algebra} and $\B$ is a \textit{$G$-$\cs$-bundle}. An action on a Hilbert $B$-module $E$ is \textit{compatible with} the $\cs$-algebra action $G \acts B$ if for each $g \in G$, $b \in B_{s(g)}$ and $e \in E_{s(g)}$, we have $g \bcdot (e \bcdot b) = (g \bcdot e) \bcdot (g \bcdot b)$, and for each $g \in G$ and $e, f \in E_{s(g)}$, we have $g \bcdot \langle e, f \rangle  = \langle g \bcdot e , g \bcdot f \rangle$. We then call $E$ a \textit{$G$-Hilbert $B$-module} and $\E$ a \textit{$G$-Hilbert $\B$-bundle}. Given $G$-Hilbert $B$-modules $E$ and $F$, conjugation by the actions of an element $g$ defines a map $T \mapsto g \bcdot T \colon \cal L(E_{s(g)},F_{s(g)}) \to \cal L(E_{r(g)},F_{r(g)})$. This determines a Banach space action $G \acts \cal K(\E,\F)$ by Lemma \ref{continuity of action map}.

Let $G$ be a Hausdorff \'etale groupoid and let $A$ be a $G$-$\cs$-algebra. The (full) crossed product $G \ltimes A$ is the completion of the $*$-algebra $\Gamma_c(G,s^*\A)$\footnote{Other authors may choose to instead use $\Gamma_c(G, r^*\A)$.} under the maximal $\cs$-norm. When $E$ is a $G$-Hilbert $A$-module, the space $\Gamma_c(G, s^*\E)$ carries a $\Gamma_c(G, s^* \A)$-valued inner product and a right action of $\Gamma_c(G,s^*\A)$. Completing with respect to this inner product and the $G \ltimes A$ norm, we obtain the crossed product Hilbert $G \ltimes A$-module $G \ltimes E$. We record the convolution and involution on $\Gamma_c(G, s^* \A)$ and the inner product and module structure on $\Gamma_c(G, s^* \E)$; for $a, b \in \Gamma_c(G, s^* \A)$, $\xi, \eta \in \Gamma_c(G, s^* \E)$ and $g \in G$ we have
\begin{align*}
a * b \colon g & \mapsto \sum_{h \in G^{s(g)}} (h \bcdot \xi(gh)) \eta(h^{-1}), & a^* \colon g & \mapsto g^{-1} \bcdot {a(g^{-1})}^*, \\
\langle \xi, \eta \rangle \colon g & \mapsto \sum_{h \in G_{r(g)}} \langle g^{-1} \bcdot \xi(h) , \eta(hg) \rangle, & \xi \bcdot a  \colon g & \mapsto \sum_{h \in G^{s(g)}} (h \bcdot \xi(gh)) \bcdot a(h^{-1}).
\end{align*}

The following lemma is a vital tool to construct representations of $\Gamma_c(G, s^* \A)$ and thus $G \ltimes A$. This is based on an important special case of \cite[Definition 5.1, Corollary 6.2]{BHM18}, see \cite[Lemma 1.83]{Miller22} for details.

\begin{lem}\label{technical lemma for crossed product}
Suppose that $G$ is a Hausdorff \'etale groupoid, $A$ is a $G$-$\Cast$-algebra, $D$ is a $\Cast$-algebra and $F$ is a Hilbert $D$-module with a dense subspace $F_0$. Suppose that $L \colon \Gamma_c(G, s^* \mathcal A) \times F_0 \to F_0$ is a bilinear map such that:
\begin{itemize}
\item the image of $L$ has dense span,
\item for $a,b \in \Gamma_c(G, s^* \A)$ and $f \in F_0$ we have $L(a,L(b,f)) = L(ab,f)$, and
\item for $a \in \Gamma_c(G, s^* \A)$ and $e, f \in F_0$ we have $\langle L(a,e), f \rangle = \langle e, L(a^*,f) \rangle$.
\end{itemize}
Then there is a unique non-degenerate $*$-representation $\pi \colon \Gamma_c(G, s^* \mathcal A) \to \mathcal L(F)$ such that $\pi(a) f = L(a,f)$ for each $a \in \Gamma_c(G, s^* \mathcal A)$ and $f \in F_0$. 
\end{lem}

A \textit{correspondence} $E \colon A \to B$ of $\Cast$-algebras is a Hilbert $B$-module $E$ with a non-degenerate $*$-representation $\varphi \colon A \to \cal L(E)$ called the \textit{structure map} of the correspondence. We typically omit $\varphi$ and write $a \bcdot e$ instead of $\varphi(a)(e)$ given $a \in A$ and $e \in E$. We call $E$ \textit{proper} if $\varphi(A) \subseteq \cal K(E)$. If $A$ and $B$ are $C_0(X)$-algebras for some locally compact Hausdorff space $X$, then we say that $E$ is a \textit{$C_0(X)$-correspondence} if the structure map is $C_0(X)$-linear. In this case, $E$ induces correspondences $E_x \colon A_x \to B_x$ of the fibres at each $x \in X$. Furthermore, if $A$ and $B$ are $G$-$\Cast$-algebras for some Hausdorff \'etale groupoid $G$ with unit space $X$ and $E$ is a $G$-Hilbert $A$-module, then $E$ is a \textit{$G$-equivariant correspondence} if 
\[g \bcdot (a \bcdot e) = (g \bcdot a) \bcdot (g \bcdot e) \]
for each $g \in G$, $a \in A_{s(g)}$ and $e \in E_{s(g)}$.

The crossed product $G \ltimes E \colon G \ltimes A \to G \ltimes B$ of a $G$-equivariant correspondence $E \colon A \to B$ is constructed as follows. Consider the left action of $\Gamma_c(G, s^* \A)$ on $\Gamma_c(G, s^* \E)$ which for $a \in \Gamma_c(G, s^* \A)$, $\xi \in \Gamma_c(G, s^* \E)$ and $g \in G$ is given by
\[
 a \bcdot \xi \colon g \mapsto \sum_{h \in G^{s(g)}} (h \bcdot a(gh)) \bcdot \xi(h^{-1}). 
\]
This extends to a non-degenerate $*$-homomorphism $G \ltimes A \to \mathcal L(G \ltimes E)$ for example by Lemma \ref{technical lemma for crossed product}. The resulting correspondence $G \ltimes E \colon G \ltimes A \to G \ltimes B$ is proper when $E$ is, because $A$ is full in $G \ltimes A$ and the action of $A$ on $G \ltimes E$ factors through $\mathcal K(E)$, which acts by compact operators on $G \ltimes E$. Explicitly, for $e, f \in E$ included as $\iota(e), \iota(f) \in G \ltimes E$ the operator $\Theta_{e,f} \in \mathcal K(E)$ acts as $\Theta_{\iota(e), \iota(f)}$ on $G \ltimes E$.

%
%

\subsection{Groupoid equivariant KK-theory}

\begin{defn}[Equivariant Kasparov cycle]
Let $G$ be a Hausdorff \'etale groupoid and let $A$ and $B$ be $G$-$\Cast$-algebras. A \textit{$G$-equivariant Kasparov $A$-$B$ cycle} is a pair $(E,T)$ where $E \colon A \to B$ is a countably generated $\Z/ 2 \Z$-graded $G$-equivariant correspondence with structure map $\varphi$, and $T \in \cal L(E)$ is an adjointable operator of degree $1$ which:
\begin{itemize}
\item is \textit{almost self-adjoint}: $\varphi(a) (T - T^*) \in \cal K(E)$ for each $a \in A$,
\item is \textit{almost unitary}: $\varphi(a) (T^*T-1) \in \cal K(E)$ for each $a \in A$,
\item \textit{almost commutes with} $\varphi$: $[T,\varphi(a)] \in \cal K(E)$ for each $a \in A$, and
\item is \textit{almost invariant}: $\varphi_{s(g)}(a) (g^{-1} \bcdot T_{r(g)} - T_{s(g)}) \in \cal K(E_{s(g)} )$ for each $g \in G$ and $a \in A_{s(g)}$, and this determines a continuous map $s^* \A \to s^* \cal K(\E)$.
\end{itemize}
We call such operators \textit{Fredholm operators}. A \textit{unitary equivalence} $(E_0,T_0) \sim_u (E_1,T_1)$ of Kasparov cycles is a $G$-equivariant unitary operator $U \in \cal L(E_0, E_1)$ of degree $0$ intertwining the actions of $A$ and the Fredholm operators $T_i$. We write $\bb E^G (A,B)$ for the set of unitary equivalence classes of $G$-equivariant Kasparov $A$-$B$ cycles, and frequently implicitly identify unitarily equivalent Kasparov cycles. A \textit{homotopy} of Kasparov cycles is a $G$-equivariant Kasparov $A$-$C([0,1],B)$ cycle, which makes the Kasparov $A$-$B$ cycles obtained by evaluation at $0$ and $1$ \textit{homotopic}. Given a Kasparov cycle $(E,T)$, a \textit{compact perturbation} of $T$ is an operator $S \in \cal L(E)$ such that $\varphi(a) (S - T) \in \cal K(E)$ and $(S - T) \varphi(a) \in \cal K(E)$ for each $a \in A$, and in this case $(E,S)$ is homotopic to $(E,T)$. The \textit{Kasparov group} $\KK^G(A,B)$ is the set of Kasparov cycles $\bb E^G(A,G)$ up to homotopy, whose elements are called \textit{Kasparov classes}. A \textit{Kasparov product} of Kasparov cycles $(E_1,T_1) \in \bb E^G(A,B)$ and $(E_2,T_2) \in \bb E^G(B,C)$ is a Kasparov cycle $(E_1 \otimes_B E_2, T)$ such that
\begin{itemize}
\item the Fredholm operator $T$ is a \textit{$T_2$-connection for $E_1$}. This means that for each homogeneous $e_1 \in E_1$, \begin{align*}
\theta_{e_1} T_2 - & (-1)^{\deg(e_1)} T \theta_{e_1} \in \cal K(E_2, E_1 \otimes_B E_2) \\
\theta_{e_1} T_2^* - & (-1)^{\deg(e_1)}  T^* \theta_{e_1} \in \cal K(E_2, E_1 \otimes_B E_2)
\end{align*}
where $\theta_{e_1} \in \cal L(E_2,E_1 \otimes_B E_2)$ is defined at $e_2 \in E_2$ by $\theta_{e_1}(e_2) = e_1 \otimes e_2$, and
\item for each $a \in A$, we have $\varphi(a)[T_1 \otimes 1, T ] \varphi(a)^* \geq 0$ mod $\cal K(E_1 \otimes_B E_2)$\footnote{The commutator is graded so in this case $[T,T_1 \otimes 1] = (T_1 \otimes 1)T + T(T_1 \otimes 1)$.}.
\end{itemize}
\end{defn}
We write $(E_1,T_1) \hash_B (E_2,T_2)$ for the subset of $\bb E^G(A,C)$ of Kasparov products of $(E_1,T_1)$ and $(E_2,T_2)$. The chief technical result in Kasparov theory is the existence and uniqueness of the Kasparov product, see \cite[6. Le produit de Kasparov]{LeGall99} for the groupoid equivariant setting.
\begin{thm}[Existence and uniqueness of the Kasparov product]
Let $G$ be a second countable Hausdorff \'etale groupoid, let $A$, $B$ and $C$ be separable $G$-$\Cast$-algebras and let $(E_1, T_1) \in \bb E^G(A,B)$ and $ (E_2,T_2) \in \bb E^G(B,C)$ be Kasparov cycles. Then there is a Kasparov product $(E,T) \in (E_1,T_1) \hash_B (E_2,T_2)$ which is unique up to homotopy. Furthermore, this descends to a bilinear operation which we call \textit{the Kasparov product}
\[ - \otimes_B - \colon \KK^G(A,B) \times \KK^G(B,C) \to \KK^G(A,C).  \]
The Kasparov product is associative and we obtain a category $\KK^G$ of separable $G$-$\Cast$-algebras whose morphisms are the Kasparov classes with composition given by the Kasparov product.
\end{thm}
In the non-equivariant setting we have $K_0(A) \cong \KK(\bb C, A)$ for each separable $\Cast$-algebra $A$. The Kasparov product then induces for each separable $A$ and $B$ a homomorphism $\KK(A,B) \to \Hom(K_*(A), K_*(B))$, which form a K-theory functor $K_* \colon \KK \to \Ab_*$. 

The crossed products $G \ltimes -$ and $G \ltimes_r -$ give rise to \textit{descent functors} $\KK^G \to \KK$ \cite[Chapitre 7]{LeGall94}, with a Kasparov cycle $(E,T) \in \bb E^G(A,B)$ sent to $(G \ltimes E, G \ltimes T)$ and $(G \ltimes_r E, G \ltimes_r T)$ respectively. The operators $G \ltimes T$ and $G \ltimes_r T$ are given by $T \otimes 1$ through identifications $G \ltimes E \cong E \otimes_B G \ltimes B$ and $G \ltimes_r E \cong E \otimes_B G \ltimes_r B$, and are both given concretely at $\nu \in \Gamma_c(G,s^* \E) \subseteq G \ltimes E$ by the element $g \mapsto (g^{-1} \bcdot T_{r(g)}) \bcdot \nu(g)$ of $\Gamma_c(G, s^* \E)$.

\section{\'Etale correspondences of \'etale groupoids}\label{correspondence section}

\'Etale correspondences are a notion of morphism between \'etale groupoids which encompass many ways through which groupoids are related, such as Morita equivalences. We refer the reader to \cite{AKM22} for the theory relevant to $\cs$-algebras. We focus here on giving many examples of \'etale correspondences to highlight their abundance in natural situations.

\begin{defn}[\'Etale correspondence]
Let $G$ and $H$ be \'etale groupoids with unit spaces $X$ and $Y$. An \textit{\'etale correspondence} $\Omega \colon G \to H$ is a space $\Omega$ with a left $G$-action and a right $H$-action with anchor maps $\rho\colon \Omega \to X$ and $\sigma \colon  \Omega \to Y$ called the \textit{range} and \textit{source} such that:
\begin{itemize}
\item the $G$-action commutes with the $H$-action - $\Omega$ is a \textit{$G$-$H$-bispace}, and
\item the right action $\Omega \rightacts H$ is free, proper and \'etale.
\end{itemize}
We say that $\Omega \colon G \to H$ is \textit{proper} if the map $\overline{\rho} \colon \Omega/H \to G^0$ induced by $\rho$ is proper. We say that $\Omega$ is \textit{Hausdorff} if each of $\Omega$, $G$ and $H$ are Hausdorff, and \textit{second countable} if all three are second countable. If we want to highlight the correspondence $\Omega$, we may write $\sigma_\Omega$ and $\rho_\Omega$ instead of $\sigma$ and $\rho$. For $x \in X$ and $y \in Y$, we write $\Omega^x$ and $\Omega_y$ for the range and source fibres $\rho^{-1}(x)$ and $\sigma^{-1}(y)$.
\end{defn}

\begin{example}[\'Etale homomorphism]\label{etale homomorphism correspondence}
A homomorphism $\varphi \colon G \to H$ of \'etale groupoids is \textit{\'etale} if it is a local homeomorphism. Let $\varphi^0 \colon G^0 \to H^0$ be the restriction of $\varphi$ to the unit space and consider the space
\[ \Omega_\varphi := G^0 \times_{\varphi^0, H^0, r} H = \{(x,h) \in G^0 \times H \suchthat \varphi(x) = r(h) \}. \]
We define a left action $G \acts \Omega_\varphi$ by $g \bcdot (s(g),h) = (r(g), \varphi(g)h)$ for $g \in G$ and $h \in H$ with $\varphi(s(g))= r(h)$ and a right action $\Omega_\varphi \rightacts H$ by $(x,h) \bcdot h' = (x,hh')$ for $x \in G^0$ and $h, h' \in H$ with $r(h) = \varphi(x)$ and $s(h) = r(h')$. The $G$-$H$-bispace $\Omega_\varphi$ is an étale correspondence when $\varphi$ is \'etale, and moreover a proper correspondence, because $\orho \colon \Omega_\varphi/H \to G^0$ is a homeomorphism with inverse given by $x \mapsto [x,\varphi(x)]_H$. An important special case of this is given by an open inclusion $\iota \colon G \hookrightarrow H$, in which case $\Omega_{\iota} = H^{G^0}$.
\end{example}

\begin{example}[Algebraic morphisms]\label{algebraic morphism example}
Let $G$ and $H$ be \'etale groupoids. An \textit{algebraic morphism} $G \to H$ is an action $G \acts H$ that commutes with the right multiplication action $H \rightacts H$. These are studied in \cite{BunSta05, BEM12} and are also known as Zakrzewski morphisms or actors. Like \'etale homomorphisms, these generalise homomorphisms of discrete groups. Unlike \'etale homomorphisms, an algebraic morphism induces a $*$-homomorphism $C^*(G) \to M(C^*(H))$ of the groupoid $\Cast$-algebras, which lands in $C^*(H)$ when the correspondence is proper.
\end{example}

\begin{example}[Topological correspondences]\label{topological correspondence}
Let $X$ and $Y$ be locally compact Hausdorff spaces considered as groupoids with only identity arrows. A \textit{topological correspondence} $\Omega \colon X \to Y$ is a locally compact Hausdorff space $\Omega$ with a continuous map $\rho \colon \Omega \to X$ and a local homeomorphism $\sigma \colon \Omega \to Y$. The correspondence is proper if and only if $\rho \colon \Omega \to X$ is proper. As a special case, consider a compact open cover $\cal U$ of a totally disconnected space $Y$. Setting $\Omega = \bigsqcup_{U \in \cal U} U$, the canonical map $\Omega \to Y$ is a local homeomorphism, and the indexing map $\Omega \to \cal U$ is proper. We obtain a proper correspondence $\Omega \colon \cal U \to Y$.
\end{example}

\begin{example}[Action correspondences]\label{action correspondence}
Let $G$ be an \'etale groupoid and let $X$ be a locally compact Hausdorff $G$-space with anchor map $\tau \colon X \to G^0$. Then there is a correspondence $G \acts G \ltimes X \rightacts G \ltimes X$ with bispace $G \ltimes X$ called the \textit{action correspondence}. The right action is given by right multiplication in $G \ltimes X$. The left action has anchor map given by $(g,x) \mapsto r(g) \colon G \ltimes X \to G^0$ and action map given by $h \bcdot (g,x) = (hg, x)$ for $h \in G$ and $(g,x) \in G \ltimes X$ with $r(g) = s(h)$. This is a special case of an algebraic morphism. The action correspondence is proper if and only if $\tau$ is a proper map. 
\end{example}

\begin{example}[Factor maps]
Let $\varphi \colon H \to G$ be a (not necessarily \'etale) continuous homomorphism of \'etale groupoids, and suppose further that $\varphi$ is surjective, proper as a continuous map and fibrewise bijective in the sense that for each $y \in H^0$, $\varphi$ restricts to a bijection $H_y \cong G_{\varphi(y)}$. In particular, the continuous bijection $\beta \colon H \to G \times_{G^0} H^0$ which sends $h \in H$ to $(\varphi(h),s(h))$ is a homeomorphism, as closedness of $\beta$ follows from properness of $\varphi$. We may then define an action $G \acts H^0$ with action map \[G \times_{G^0} H^0 \xrightarrow{\beta^{-1}} H \xrightarrow{r} H^0.\]
The map $\beta$ becomes a groupoid isomorphism $\beta \colon H \to G \ltimes H^0$, and so we obtain a proper action correspondence from $G$ to $H$. Factor maps such as $\varphi$ are of interest (for example in \cite{Li20}) because they induce an inclusion $C^*_r(G) \subseteq C^*_r(H)$.
\end{example}

We can also build a correspondence of groupoids using only the ``right hand side'' data.
\begin{example}\label{correspondence of a free proper etale space}
Let $H$ be an \'etale groupoid and let $\Omega$ be a free, proper, \'etale right $H$-space. Then $\Omega$ is a proper correspondence from $\Omega/H$ to $H$.
\end{example}

Let $S \acts X$ and $S \acts Y$ be actions of an inverse semigroup by partial homeomorphisms on locally compact Hausdorff spaces $X$ and $Y$ and let $f \colon X \to Y$ be an $S$-equivariant continuous map. This induces a homomorphism $S \ltimes X \to S \ltimes Y$ sending $[s,x] \in S \ltimes X$ to $[s,f(x)] \in S \ltimes Y$ which is \'etale if $f$ is a local homeomorphism. On the other hand, if $f^{-1}(\dom_Y s) = \dom_X s$ for each $s \in S$, then there is an action $S \ltimes Y \acts S \ltimes X$ by left multiplication in $S$. Combining these, we may form \'etale correspondences from suitably $S$-equivariant topological correspondences

\begin{defn}
Let $S$ be an inverse semigroup and let $\Omega \colon X \to Y$ be a topological correspondence equipped with actions $S \acts \Omega$, $S \acts X$ and $S \acts Y$ by partial homeomorphisms. We say that $\Omega$ is an \textit{$S$-equivariant topological correspondence} if the range and source $\rho \colon \Omega \to X$ and $\sigma \colon \Omega \to Y$ are $S$-equivariant and for each $s \in S$, we have $\rho^{-1}(\dom_X s) = \dom_\Omega s$.
\end{defn}

\begin{example}[Inverse semigroup equivariant topological correspondence]
Let $S$ be an inverse semigroup and let $\Omega \colon X \to Y$ be an $S$-equivariant topological correspondence. Then the space 
\[ \Omega \times_{\sigma, Y, r} (S \ltimes Y) = \{(\omega,[s,y]) \in \Omega \times (S \ltimes Y) \suchthat \sigma(\omega) = s \bcdot y  \} \]
may be given the structure of an \'etale correspondence from $S \ltimes X$ to $S \ltimes Y$. The left action is given by $[s,x] \bcdot (\omega,[t,y]) = (s \bcdot \omega, [st,y])$ for $[s,x] \in S \ltimes X$, $\omega \in \Omega$ and $[t,y] \in S \ltimes Y$ with $\rho(\omega) = x$ and $\sigma(\omega) = t \bcdot y$. The right action is given by $(\omega,[s,y]) \bcdot [t,z] = (\omega, [st,z])$ for $\omega \in \Omega$ and $[s,y], [t,z] \in S \ltimes Y$ with $\sigma(\omega) = s \bcdot y$ and $y = t \bcdot z$. If $\Omega \colon X \to Y$ is proper, then so is $S \ltimes \Omega \colon S \ltimes X \to S \ltimes Y$.
\end{example}

\begin{example}
Let $S$ be an inverse semigroup with idempotent semilattice $E$. Let $E^\times$ be the set of non-zero idempotents and let $\hat E$ be the space of filters on $E$, with compact open sets $U_e = \{ \chi \in \hat E \suchthat \chi(e) = 1 \}$ for each $e \in E^\times$. There is a canonical action $S \acts E^\times$ given by $s \bcdot e = s e s^*$ for $s \in S$ and $e \in E$ with $e s^* s = e$, and a canonical action $S \acts \hat E$ given by $s \bcdot \chi \colon e \mapsto \chi(s^* e s)$ for $s \in S$, $\chi \in U_{s^* s}$ and $e \in E$. Then the topological correspondence
\[ E^\times \leftarrow \bigsqcup_{e \in E^\times} U_e \to \hat E \]
is an $S$-equivariant topological correspondence, with $s \bcdot (e,\chi) \defeq (s \bcdot e, s \bcdot \chi)$ for $s \in S$, $e \in \dom_{E^\times} s$ and $\chi \in U_e$. This induces a correspondence
\[ S \ltimes E^\times \acts \bigsqcup_{e \in E^\times} \left\{ [s,\chi] \in S \ltimes \hat E \suchthat s \bcdot \chi \in U_e \right\} \rightacts S \ltimes \hat E. \]
\end{example}

To compose \'etale correspondences $\Omega \colon G \to H$ and $\Lambda \colon H \to K$ we take the pullback $\Omega \times_H \Lambda$ of $\Omega$ and $\Lambda$ over $H$. 

\begin{defn}[Composition of \'etale correspondences]
Let $\Omega \colon G \to H$ and $\Lambda \colon H \to K$ be \'etale correspondences. The \textit{composition} $\Lambda \circ \Omega \colon G \to K$ is given by the pullback
\[ \Omega \times_H \Lambda = (\Omega \times_{H^0} \Lambda)/H \]
which is naturally a $G$-$K$-bispace. The left $G$-action is induced by the action on $\Omega$ and the right $K$-action is induced by the action on $\Lambda$. Explicitly, the actions are given by the following formulae for $\omega \in \Omega$, $\lambda \in \Lambda$, $g \in G$ and $k \in K$ with $s(g) = \rho(\omega)$, $\sigma(\omega) = \rho(\lambda)$ and $\sigma(\lambda) = r(k)$.
\begin{align*}
G \acts \Omega & \times_H \Lambda & g \bcdot  [\omega,\lambda]_H & \phantom{ {} \bcdot k } := [g \bcdot \omega, \lambda]_H, \\
  \Omega & \times_H \Lambda \rightacts K & [\omega,\lambda]_H & \bcdot k := [\omega, \lambda \bcdot k]_H.
\end{align*}
\end{defn}

\'Etale groupoids and \'etale correspondences naturally form a bicategory in which we may further consider 2-arrows between \'etale correspondences \cite[Section 6]{AKM22}, as composition is associative up to canonical isomorphisms. For our purposes, it is enough to ignore 2-categorical subtleties and instead consider the category of \'etale groupoids whose morphisms are isomorphism classes of \'etale correspondences.

\begin{example}[Morita equivalences]
Let $G$ and $H$ be \'etale groupoids. A \textit{Morita equivalence} from $G$ to $H$ is a $G$-$H$-bispace $\Omega$ with anchor maps $\rho \colon \Omega \to G^0$ and $\sigma \colon \Omega \to H^0$ such that
\begin{itemize}
\item the left action $G \acts \Omega$ and the right action $\Omega \rightacts H$ are free and proper, and
\item the maps $\overline{\rho} \colon \Omega/H \to G^0$ and $\overline{\sigma} \colon G \backslash \Omega \to H^0$ induced by the anchor maps are homeomorphisms.
\end{itemize}
The second condition implies that the actions are \'etale, so this is a (proper) correspondence. Morita equivalences are exactly the invertible \'etale correspondences \cite[Theorem 2.30]{Albandik15}. 
\end{example}

\begin{example}[Étale Hilsum--Skandalis morphisms]\label{HS morphism}
A \textit{Hilsum--Skandalis morphism} $\Lambda \colon G \to H$ of étale groupoids is a $G$-$H$-bispace $\Lambda$ with anchor maps $\rho \colon \Lambda \to G^0$ and $\sigma \colon \Lambda \to H^0$ such that
\begin{itemize}
\item the right action $\Lambda \rightacts H$ is free and proper, and
\item the map $\orho \colon \Lambda/H \to G^0$ induced by the left anchor map is a homeomorphism.
\end{itemize}
A Hilsum--Skandalis morphism of étale groupoids which is étale in the sense that $\sigma \colon \Lambda \to H^0$ is étale is an étale correspondence. Conversely, any étale correspondence $\Omega \colon G \to H$ decomposes into the action correspondence $G \to G \ltimes \Omega/H$ and the étale Hilsum--Skandalis morphism $G \ltimes \Omega/H \acts \Omega \rightacts H$. Le Gall refers to Hilsum--Skandalis morphisms as regular graphs, for which his generalised morphisms provide an alternative perspective \cite[Définition 2.3]{LeGall99}.
\end{example}

When working with an \'etale groupoid $G$, the open bisections $U \subseteq G$ are particularly important, and we may refer to them as \textit{slices}. Given an \'etale correspondence $\Omega \colon G \to H$, a \textit{slice} $U \subseteq \Omega$ is an open set such that $q \colon \Omega \to \Omega/H$ and $\sigma \colon \Omega \to H^0$ are injective on $U$. If we have another \'etale correspondence $\Lambda \colon H \to K$, the composition $V \circ U \subseteq \Omega \times_H \Lambda$ of slices $U \subseteq \Omega$ and $V \subseteq \Lambda$ is the slice
\[ V \circ U := \{ [u,v]_H \in \Omega \times_H \Lambda \suchthat (u,v) \in U \times_{H^0} V \}. \]

Given an \'etale correspondence $\Omega \colon G \to H$ there is a correspondence 
\[C^*(\Omega) \colon C^*(G) \to C^*(H)\]
of the full $\Cast$-algebras\cite[Section 7]{AKM22} obtained as a completion of $C_c(\Omega)$\footnote{When $X$ is not Hausdorff, the notation $C_c(X)$ is an abuse of notation which refers to the set of (not necessarily continuous) functions $X \to \bb C$ which may be written as a finite sum of functions which vanish outside of an open Hausdorff subspace $U \subseteq X$ and are continuous and compactly supported on $U$.}. This construction respects composition of correspondences. The correspondence $C^*(\Omega)$ is proper when $\Omega$ is proper, and thus induces a map in K-theory.

\section{Induced Banach spaces}\label{induced banach space section}

In this section we will describe how an \'etale correspondence $\Omega \colon G \to H$ may be used to induce a $G$-$\Cast$-algebra from an $H$-$\Cast$-algebra. In order to do this at the level of $\KK$-theory, we also cover induced Hilbert modules and spaces of compact operators through the lens of $H$ and $G$-Banach spaces.

\begin{defn}[Induced Banach space]
Let $\Omega \colon G \to H$ be a Hausdorff \'etale correspondence and let $A$ be an $H$-Banach space. We define the induced Banach space $\Ind_\Omega A$ to be the space of sections $\xi \in \Gamma_b(\Omega,\sigma^* \A)$ such that:
\begin{itemize}
\item for each $\omega \in \Omega$ and $h \in H^{\sigma(\omega)}$, we have $\xi(\omega \bcdot h) = h^{-1} \bcdot \xi(\omega)$, and
\item the map $[\omega]_H \mapsto \lVert \xi(\omega) \rVert \colon \Omega/H \to \mathbb R_{\geq 0}$ vanishes at infinity.
\end{itemize}
This is a Banach space with the norm $\lVert \xi \rVert = \sup_{\omega \in \Omega} \lVert \xi(\omega) \rVert$. 
\end{defn}

\begin{prop}
Let $\Omega \colon G \to H$ be a Hausdorff \'etale correspondence and let $A$ be an $H$-Banach space. Then $\Ind_\Omega A$ may be given the structure of a $C_0(\Omega/H)$-Banach space and a $C_0(G^0)$-Banach space with the following structure maps. For $\gamma \in C_0(\Omega/H)$, $\delta \in C_0(G^0)$ and $\xi \in \Ind_\Omega A$ we define $(\gamma \bcdot \xi)(\omega) = \gamma([\omega]_H)\xi(\omega)$ and $(\delta \bcdot \xi)(\omega) = \delta(\rho(\omega))\xi(\omega)$ for each $\omega \in \Omega$.
\begin{proof}
We need only check the non-degeneracy of the actions of $C_0(\Omega/H)$ and $C_0(G^0)$ on $\Ind_\Omega A$. For each $\xi \in \Ind_\Omega A$ and $\epsilon > 0$ there is a compact set $K \subseteq \Omega/H$ such that $\lVert \xi(\omega) \rVert < \epsilon$ for $\omega \in \Omega$ with $[\omega]_H$ outside $K$. There are $[0,1]$-valued functions $\gamma \in C_c(\Omega/H)$ and $\delta \in C_c(G^0)$ such that $\gamma = 1$ on $K$ and $\delta = 1$ on $\overline{\rho}(K)$. By construction, we have $\lVert \gamma \bcdot \xi - \xi \rVert < \epsilon$ and $\lVert \delta \bcdot \xi - \xi \rVert < \epsilon$.
\end{proof}
\end{prop}

We may construct an abundance of elements of induced Banach spaces by extending sections with compact support on some slice $U \subseteq \Omega$ in an $H$-equivariant fashion.

\begin{prop}\label{equivariant extensions}
Let $\Omega \colon G \to H$ be a Hausdorff \'etale correspondence and let $A$ be an $H$-Banach space. Suppose that $\eta \in \Gamma_c(\Omega,\sigma^* \A)$ is compactly supported on a slice $U \subseteq \Omega$. Then the $H$-equivariant extension $\tilde \eta \colon \Omega \to \sigma^* \A$ of $\eta$ defined by
\[ \tilde \eta \colon \omega \mapsto \sum_{h \in H^{\sigma(\omega)}} h \bcdot \eta(\omega \bcdot h) \]
is an element of $\Ind_\Omega A$ with $\lVert \tilde \eta \rVert = \lVert \eta \rVert_\infty$. If $\cal U$ is an open cover of $\Omega$ by slices then elements of this form for $U \in \cal U$ have dense span in $\Ind_\Omega A$.
\begin{proof}
Consider the range map $r \colon \Omega \rtimes H \to \Omega$ and the continuous section 
\[ (\omega,h) \mapsto h \bcdot \eta(\omega \bcdot h) \colon \Omega \rtimes H \to r^* \sigma^* \A. \]
Because $\eta$ is compactly supported and $\Omega \rtimes H$ is proper, this section is properly supported with respect to $r \colon \Omega \rtimes H \to \Omega$. Lemma \ref{propersupportlemma} therefore justifies that $\tilde \eta$ is well-defined and continuous. It is $H$-equivariant and its support has compact image in $\Omega/H$ by construction, so defines an element of $\Ind_\Omega A$. As $\eta$ is supported on a slice, the sum at each $\omega \in \Omega$ can have at most one non-zero summand, and therefore $\lVert \tilde \eta \rVert = \lVert \eta \rVert_\infty$.

The subspace $C_c(\Omega/H) \Ind_\Omega A$ is dense in $\Ind_\Omega A$ and $C_c(\Omega/H)$ is spanned by $C_c(q(U))$ across all $U \in \cal U$, where $q \colon \Omega \to \Omega/H$ is the quotient map. Any element of $C_c(q(U)) \Ind_\Omega A$ is the $H$-equivariant extension of its restriction to $U$, and therefore $H$-equivariant extensions have dense span in $\Ind_\Omega A$.
\end{proof}
\end{prop}

\begin{prop}\label{fibres of induced Banach spaces}
Let $\Omega \colon G \to H$ be a Hausdorff \'etale correspondence, let $A$ be an $H$-Banach space and let $x \in G^0$. Then the fibre $(\Ind_\Omega A)_x$ at $x \in G^0$ of $\Ind_\Omega A$ as a $C_0(G^0)$-Banach space is the space of sections $\xi \in \Gamma_b(\Omega^x, \sigma^* \A)$ such that:
\begin{itemize}
\item for each $\omega \in \Omega^x$ and $h \in H^{\sigma(\omega)}$, we have $\xi(\omega \bcdot h) = h^{-1} \bcdot \xi(\omega)$,
\item and the map $[\omega]_H \mapsto \lVert \xi(\omega) \rVert \colon \Omega^x/H \to \mathbb R_{\geq 0}$ vanishes at infinity.
\end{itemize}
The fibre map $\Ind_\Omega A \to (\Ind_\Omega A)_x$ is given by restriction from $\Omega$ to $\Omega^x$. For each $\omega \in \Omega$, the fibre $(\Ind_\Omega A)_{[\omega]_H}$ at $[\omega]_H \in \Omega/H$ of $\Ind_\Omega A$ as a $C_0(\Omega/H)$-Banach space is isomorphic to $A_{\sigma(\omega)}$, with the fibre map $\Ind_\Omega A \to A_{\sigma(\omega)}$ given by evaluation at $\omega$.
\begin{proof} For each $x \in G^0$, the restriction map $\Ind_\Omega A \to \Gamma_b(\Omega^x, \sigma^* \A)$ vanishes on $C_0(G^0 \setminus \{x\}) \Ind_\Omega A$. Conversely, if $\xi \in \Ind_\Omega A$ vanishes on $\Omega^x$, then for each $\epsilon > 0$, let $K = \{ \omega \in \Omega \suchthat \lVert \xi(\omega) \rVert \geq \epsilon\}$. The set $\rho(K) \subseteq G^0$ is compact and does not contain $x$. We may therefore find a $[0,1]$-valued function $\gamma \in C_c(G^0 \setminus \{x\})$ such that $\gamma = 1$ on $\rho(K)$. By construction, $\lVert \gamma \bcdot \xi - \xi \rVert < \epsilon$ and therefore $\overline{C_0(G^0 \setminus \{x\}) \Ind_\Omega A}$ is the kernel of the restriction map $\Ind_\Omega A \to \Gamma_b(\Omega^x, \sigma^* \A)$. The fibre $(\Ind_\Omega A)_x$ may therefore be identified with the image of the restriction map. For each element of $\Ind_\Omega A$, its restriction to $\Omega^x$ is $H$-equivariant and vanishes at infinity with respect to $\Omega^x/H$. For each slice $U \subseteq \Omega$, the $H$-equivariant extension of an element $\eta \in \Gamma_c(U, \sigma^* \A)$ restricts to the $H$-equivariant extension of the restriction $\eta_x \in \Gamma_c(U \cap \Omega^x, \sigma^* \A)$ of $\eta$ to $\Omega^x$. These elements have dense span in the Banach space of $H$-equivariant sections in $\Gamma_b(\Omega^x, \sigma^* \A)$ which vanish at infinity with respect to $\Omega^x/H$. This space is therefore the image of the restriction map and so identified with $(\Ind_\Omega A)_x$.

The description of the fibres of $\Ind_\Omega A$ over $\Omega/H$ are a special case of the description over $G^0$, because $\Omega$ is also an \'etale correspondence from $\Omega/H$ to $H$. In this setting, the $H$-equivariant sections in $\Gamma_b(\Omega^{[\omega]_H}, \sigma^* \A)$ which vanish at infinity with respect to $\Omega^{[\omega]_H}/H = \{ [\omega]_H \}$ may be identified through evaluation at $\omega$ with $A_{\sigma(\omega)}$. 
\end{proof}
\end{prop}
As the fibre maps are given by restrictions and evaluations, then $\Ind_\Omega A$ is locally $C_0(G^0)$-convex as a $C_0(G^0)$-Banach space and locally $C_0(\Omega/H)$-convex as a $C_0(\Omega/H)$-Banach space. This provides us with a convenient criterion for density in $\Ind_\Omega A$.
\begin{prop}\label{density in induced Banach space}
Let $\Omega \colon G \to H$ be a Hausdorff \'etale correspondence and let $A$ be an $H$-Banach space. Let $\Gamma \subseteq \Ind_\Omega A$ be a vector subspace closed under the action of $C_0(\Omega/H)$ such that $\{ \xi(\omega) \suchthat \xi \in \Gamma \}$ is dense in $A_{\sigma(\omega)}$ for each $\omega \in \Omega$. Then $\Gamma$ is dense in $\Ind_\Omega A$.
\begin{proof}
Local $C_0(\Omega/H)$-convexity of $\Ind_\Omega A$ allows us to view it as a space of $c_0$ sections, and therefore we may apply the criterion for density in Proposition \ref{density in Banach spaces}. The description of the fibres of $\Ind_\Omega A$ over $\Omega/H$ in Proposition \ref{fibres of induced Banach spaces} completes the proof.
\end{proof}
\end{prop}
As the fibres of $\Ind_\Omega A$ over $G^0$ are themselves section spaces, any map into the associated bundle $\Ind_\Omega \A$ over $G^0$ may again be evaluated at various points in $\Omega$. Continuity of these maps may be then be understood in terms of an associated bivariate function which lands in $\A$. This is typically the most useful way to get a handle on the topology of $\Ind_\Omega \A$.
\begin{lem}[Continuity of maps into induced bundles]\label{maps into induced bundle}
Let $\Omega \colon G \to H$ be a Hausdorff \'etale correspondence, let $A$ be an $H$-Banach space and consider the induced bundle $p \colon \Ind_\Omega \A \to G^0$. A map $\xi \colon Z \to \Ind_\Omega \A$ from a locally compact Hausdorff space $Z$ is continuous if and only if
\begin{itemize}
\item the composition $g = p \circ \xi \colon Z \to G^0$ is continuous, and
\item the map $\tilde \xi \colon (z, \omega) \mapsto \gamma(z) \xi(z)(\omega) \colon Z \times_{G^0} \Omega \to \A$ is $c_0$ with respect to $Z \times_{G^0} \Omega/H$ for each $\gamma \in C_c(Z)$.
\end{itemize}
\begin{proof}
We first observe that the map 
\begin{equation}\label{induced bundle evaluation continuous map}
(\eta,\omega) \mapsto \eta(\omega) \colon \Ind_\Omega \A \times_{G^0} \Omega \to \A
\end{equation} 
is continuous. Indeed, for a convergent net $(\eta_i, \omega_i) \to (\eta_0, \omega_0)$ in $\Ind_\Omega \A \times_{G^0} \Omega$, we may pick $\nu \in \Ind_\Omega A$ with $\nu_{\rho(\omega_0)} = \eta_0$. We then have convergent nets $\eta_i \to \eta_0$ and $\nu_{\rho(\omega_i)} \to \eta_0$ and so $\eta_i(\omega_i) - \nu(\omega_i) \to 0$. Since $\nu(\omega_i) \to \eta_0(\omega_0)$, we have $\eta_i(\omega_i) \to \eta_0(\omega_0)$. 

Suppose that $\xi \colon Z \to \Ind_\Omega \A$ is continuous. Then $g = p \circ \xi$ is the composition of two continuous functions. Let $\pi_Z$ and $\pi_\Omega$ be the coordinate projections on $Z \times_{G^0} \Omega$. For any $\mu \in \Gamma_0(Z, g^* \Ind_\Omega \A)$, we obtain a continuous section $(z,\omega) \mapsto \mu(z)(\omega) \colon Z \times_{G^0} \Omega \to \pi_\Omega^* \sigma^* \A$ by continuity of (\ref{induced bundle evaluation continuous map}). Each $\mu$ may be approximated by elements of the form $\gamma \eta \colon z \mapsto \gamma(z) \eta_{g(z)}$ for $\gamma \in C_c(Z)$ and $\eta \in \Ind_\Omega A$ by Proposition \ref{density in Banach spaces}. These elements are mapped to bounded sections $Z \times_{G^0} \Omega \to \pi_\Omega^* \sigma^* \A$ which vanish at infinity with respect to $Z \times_{G^0} \Omega/H$. By approximation, the section $(z,\omega) \mapsto \mu(z)(\omega) \colon Z \times_{G^0} \Omega \to \pi^*_\Omega \sigma^* \A$ vanishes at infinity for each $\mu \in \Gamma_0(Z, g^* \Ind_\Omega \A)$. Applying this to $\mu = \gamma \xi \colon z \mapsto \gamma(z) \xi(z)$, we may conclude that $\tilde \xi$ is $c_0$ with respect to $Z \times_{G^0} \Omega/H$. 

Now suppose that $g = p \circ \xi$ is continuous and that $\tilde \xi$ is continuous for each $\gamma \in C_c(Z)$. We will show that $\gamma \xi$ is continuous when $\lVert \gamma \rVert_\infty \leq 1$ and therefore $\xi$ must be continuous. Let $z_i \to z_0$. Let $x_i = g(z_i)$, $x_0 = g(z_0)$ and let $\eta \in \Ind_\Omega A$ with $\eta_{x_0} = \xi(z_0)$. We set $u_i = \gamma(z_i)\eta_{x_i}$, which converge to $\gamma(z_0) \xi(z_0)$. For each $\epsilon > 0$, the set
\begin{align*}
& \left\{ z \in Z \suchthat  \lVert \gamma(z)( \xi(z) - \eta_{g(z)})\rVert \geq \epsilon \right\} \\
& = \left\{ z \in Z \suchthat \sup_{\omega \in \Omega^{g(z)}} \lVert \gamma(z)( \xi(z)(\omega) - \eta(\omega)) \rVert \geq \epsilon \right\} \\
& = \bigcap_{\delta \in (0,\epsilon)} \pi_Z \left( \left\{ (z,[\omega]_H) \in Z \times_{G^0} \Omega/H \suchthat \lVert \gamma(z)( \xi(z)(\omega) - \eta(\omega))\rVert \geq \delta \right\} \right)
\end{align*}
is the intersection of continuous images of compact sets and so closed. For large $i$ we therefore have $\lVert \gamma(z_i) \xi(z_i) - u_i \rVert < \epsilon$ and so $\gamma(z_i) \xi(z_i) - u_i \to 0$. We may conclude that $\gamma(z_i) \xi(z_i) \to \gamma(z_0)\xi(z_0)$.
\end{proof}
\end{lem}

We may use this lemma to construct an action of $G$ on $\Ind_\Omega A$, so that $\Ind_\Omega$ actually induces $G$-Banach spaces from $H$-Banach spaces.

\begin{prop}
Let $\Omega \colon G \to H$ be a Hausdorff \'etale correspondence and let $A$ be an $H$-Banach space. Then the following is a (continuous) action of $G$ on the Banach bundle $\Ind_\Omega \A \to G^0$. For $g \in G$ and $\xi \in (\Ind_\Omega A)_{s(g)}$, the element $g \bcdot \xi \in (\Ind_\Omega A)_{r(g)}$ is given by the section
\[ g \bcdot \xi \colon \omega \mapsto \xi(g^{-1} \bcdot \omega) \colon \Omega^{r(g)} \to \sigma^* \A. \] 
\begin{proof}
By Lemma \ref{continuity of action map} we may check that the action map $(g,\xi) \mapsto g \bcdot \xi \colon s^* \Ind_\Omega \A \to r^* \Ind_\Omega \A$ is continuous by checking that for each $\eta \in \Ind_\Omega A$, the map $g \mapsto g \bcdot \eta_{s(g)} \colon G \to \Ind_\Omega \A$ is continuous. By Lemma \ref{maps into induced bundle} this reduces to checking that for each $\gamma \in C_c(G)$, the map $G \times_{r,G^0,\rho} \Omega \to \A$ given by 
\[ (g, \omega) \mapsto \gamma(g)\eta(g^{-1} \bcdot \omega) \]
is $c_0$ with respect to $G \times_{G^0} \Omega/H$. Let $\epsilon > 0$ and consider the compact set $K = \{ [\omega]_H \in \Omega/H \suchthat \lVert \gamma \rVert_{\infty} \lVert \eta(\omega) \rVert \geq \epsilon \}$. The above map is $\epsilon$-bounded away from the compact set $\supp \gamma \times_{G^0} (\supp \gamma)^{-1} \bcdot K \subseteq G \times_{G^0} \Omega/H$.
\end{proof}
\end{prop}
Now that we have all the necessary general tools for working with induced Banach spaces, we turn our attention to $\Cast$-algebras.
\begin{prop}
Let $\Omega \colon G \to H$ be a Hausdorff \'etale correspondence and let $A$ be an $H$-$\Cast$-algebra. Then the induced $G$-Banach space $\Ind_\Omega A$ is a $G$-$\Cast$-algebra.
\begin{proof}
It is clear through pointwise computations over $\Omega$ that $C_0(G^0)$ acts by central multipliers on $\Ind_\Omega A$, turning it into a $C_0(G^0)$-algebra. By checking pointwise over $\Omega^{s(g)}$, each $g \in G$ induces a $*$-homomorphism $(\Ind_\Omega A)_{s(g)} \to (\Ind_\Omega A)_{r(g)}$ of fibres.
\end{proof}
\end{prop}

\begin{prop}
Let $\Omega \colon G \to H$ be a Hausdorff \'etale correspondence, let $A$ be an $H$-$\Cast$-algebra and let $E$ be an $H$-Hilbert $A$-module. Then the induced $G$-Banach space $\Ind_\Omega E$ is a Hilbert $\Ind_\Omega A$-module, with operations given pointwise over $\Omega$. The Banach space action $G \acts \Ind_\Omega E$ is compatible with the $\Cast$-algebra action $G \acts \Ind_\Omega A$.
\begin{proof}
The pointwise operations on $\Ind_\Omega E$ are defined for $\xi, \eta \in \Ind_\Omega E$ and $\nu \in \Ind_\Omega A$ by $\langle \xi, \eta \rangle \colon \omega \mapsto \langle \xi(\omega), \eta(\omega) \rangle \colon \Omega \to \sigma^* \A$ and $\xi \bcdot \nu \colon \omega \mapsto \xi(\omega) \bcdot \nu(\omega) \in \sigma^* \E$ for each $\omega \in \Omega$. These are bounded, $H$-equivariant and vanish at infinity with respect to $\Omega/H$, and continuous by Remark \ref{bilinear continuity remark}. Thus $\langle \xi, \eta \rangle \in \Ind_\Omega A$ and $\xi \bcdot \nu \in \Ind_\Omega E$ and we obtain the structure of an inner product $\Ind_\Omega A$-module on $\Ind_\Omega E$. This induces the same norm that $\Ind_\Omega E$ came with, so turns $\Ind_\Omega E$ into a Hilbert $\Ind_\Omega A$-module.

Implicit in the statement of the proposition is that the $C_0(G^0)$-Banach space structure on $\Ind_\Omega E$ agrees with the $C_0(G^0)$-Banach space structure it inherits as a Hilbert $\Ind_\Omega A$-module, which is readily verified. Compatibility of $G \acts \Ind_\Omega E$ with $G \acts \Ind_\Omega A$ follows from the actions only permuting the input which the operations are taken pointwise over.
\end{proof}
\end{prop}

\begin{prop}[Fibres of an adjointable operator on induced modules]
Let $\Omega \colon G \to H$ be a Hausdorff \'etale correspondence, let $A$ be an $H$-$\Cast$-algebra, let $E$ and $F$ be $H$-Hilbert $A$-modules and let $T \in \cal L(\Ind_\Omega E, \Ind_\Omega F)$. For each $\omega \in \Omega$, there is a unique operator $T_\omega \in \cal L(E_{\sigma(\omega)},F_{\sigma(\omega)})$ such that $T_\omega \xi(\omega) = (T \xi)(\omega)$ for each $\xi \in \Ind_\Omega E$. 
\begin{proof}
Given $e \in E_{\sigma(\omega)}$ we define $T_\omega (e)$ to be $(T \xi)(\omega)$ for any $\xi \in \Ind_\Omega E$ with $\xi(\omega) = e$. This is independent of the choice of $\xi$ because for any $\eta \in \Ind_\Omega F$, we have $\langle (T \xi)(\omega), \eta(\omega) \rangle = \langle T \xi, \eta \rangle(\omega) = \langle \xi, T^* \eta \rangle(\omega) = \langle e , (T^* \eta) (\omega) \rangle$. We obtain an operator $T_\omega \colon E_{\sigma(\omega)} \to F_{\sigma(\omega)}$ and applying the same argument to $T^*$ we obtain an operator $T^*_\omega \colon F_{\sigma(\omega)} \to E_{\sigma(\omega)}$ which is adjoint to $T_\omega$. The formula $T_\omega \xi(\omega) = (T \xi)(\omega)$ for each $\xi \in \Ind_\Omega E$ is satisfied by construction and determines $T_\omega$ uniquely.
\end{proof}
\end{prop}
We call $T_\omega$ the \textit{fibre} of $T \in \cal L(\Ind_\Omega E, \Ind_\Omega F)$ at $\omega \in \Omega$. The operator $T$ is determined by its fibres and moreover $\lVert T \rVert = \sup_{\omega \in \Omega} \lVert T_\omega \rVert$. For $\xi \in \Ind_\Omega E$ and $\eta \in \Ind_\Omega F$, the fibre $(\Theta_{\eta,\xi})_\omega$ is given by $\Theta_{\eta(\omega),\xi(\omega)}$.

\begin{prop}[Description of compact operators on induced modules]\label{compactoperatorsinducedmodule}
Let $\Omega \colon G \to H$ be a Hausdorff \'etale correspondence, let $A$ be an $H$-$\Cast$-algebra and let $E$ and $F$ be $H$-Hilbert $A$-modules.
For each $T \in \cal K(\Ind_\Omega E, \Ind_\Omega F)$, the section $\omega \mapsto T_\omega \colon \Omega \to \sigma^* \cal K(\E,\F)$ defines an element of $\Ind_\Omega \cal K(E,F)$. Furthermore, this describes an isometric isomorphism $\cal K(\Ind_\Omega E, \Ind_\Omega F) \cong \Ind_\Omega \cal K(E,F)$. 
\begin{proof}
It is straightforward to see that the section $\omega \mapsto T_\omega$ is bounded and $H$-equivariant. It suffices to prove that it is continuous and vanishes at infinity with respect to $\Omega/H$ when $T = \Theta_{\eta,\xi}$ for $\xi \in \Ind_\Omega E$ and $\eta \in \Ind_\Omega F$, as elements of this form have dense span in $\cal K(\Ind_\Omega E, \Ind_\Omega F)$. The map $[\omega]_H \mapsto \lVert \Theta_{\eta(\omega),\xi(\omega)} \rVert \colon \Omega/H \to \bb R_{\geq 0}$ vanishes at infinity because $\lVert \Theta_{\eta(\omega),\xi(\omega)} \rVert \leq \lVert \eta(\omega) \rVert \lVert \xi(\omega) \rVert$ for each $\omega \in \Omega$. The continuity of $\omega \mapsto \Theta_{\eta(\omega),\xi(\omega)} \colon \Omega \to \sigma^* \cal K(\E,\F)$ follows from the continuity of $(f,e) \mapsto \Theta_{f,e} \colon \F \times_{H^0} \E \to \cal K(\E,\F)$.

This defines an isometry $\varphi \colon \cal K(\Ind_\Omega E, \Ind_\Omega F) \to \Ind_\Omega \cal K(E,F)$. The image of $\varphi$ is closed under the action of $C_0(\Omega/H)$. For each $\omega \in \Omega$, the image in $\cal K(E_{\sigma(\omega)},F_{\sigma(\omega)})$ contains $\Theta_{\eta(\omega),\xi(\omega)}$ for each $\xi \in \Ind_\Omega E$ and $\eta \in \Ind_\Omega F$. These elements have dense span, and so $\varphi$ is surjective by Proposition \ref{density in induced Banach space}.
\end{proof}
\end{prop}

\begin{prop}[Description of adjointable operators on induced modules]\label{adjointableoperatorsinducedmodule}
Let $\Omega \colon G \to H$ be a Hausdorff \'etale correspondence, let $A$ be an $H$-$\Cast$-algebra and let $E$ and $F$ be $H$-Hilbert $A$-modules. For each $T \in \cal L(\Ind_\Omega E, \Ind_\Omega F)$, the section $\omega \mapsto T_\omega \colon \Omega \to \cal L(\sigma^* \E, \sigma^* \F)$ is  bounded, strictly continuous and $H$-equivariant. Furthermore, this describes an isometry $\cal L(\Ind_\Omega E, \Ind_\Omega F) \hookrightarrow \Gamma_b(\Omega, \cal L(\sigma^* \E,\sigma^* \F))$ whose image is the $H$-equivariant sections.
\begin{proof}
It is straightforward to see that the section $\omega \mapsto T_\omega$ is bounded and $H$-equivariant. To check that $(T_\omega)_{\omega \in \Omega} \colon \sigma^* \E \to \sigma^* \E$ is continuous, we may check that for each slice $U \subseteq \Omega$ and each $\eta \in \Gamma_c(U, \sigma^* \E)$, the section $u \mapsto T_u \eta(u)\colon U \to \sigma^* \E$ is continuous. This section is the restriction of $T(\tilde \eta) \in \Ind_\Omega E$ to $U$, where $\tilde \eta \in \Ind_\Omega E$ of $\eta$ is the $H$-equivariant extension of $\eta$. It is therefore continuous, and so $(T_\omega)_{\omega \in \Omega} \colon \sigma^* \E \to \sigma^* \E$ is continuous. The adjoint $(T^*_\omega)_{\omega \in \Omega} \colon \sigma^* \E \to \sigma^* \E$ is similarly continuous and therefore $\omega \mapsto T_\omega$ is strictly continuous.

To reconstruct an adjointable operator from a bounded, strictly continuous and $H$-equivariant section $\omega \mapsto T_\omega \colon \Omega \to  \cal L(\sigma^* \E,\sigma^* \F)$, we define $T(\xi) \colon \omega \mapsto T_\omega(\xi(\omega))$ for $\xi \in \Ind_\Omega E$. This is continuous, $H$-equivariant and vanishes at infinity with respect to $\Omega/H$ by construction, so defines an element of $\Ind_\Omega F$. We obtain an operator $T \colon \Ind_\Omega E \to \Ind_\Omega F$ and applying the same argument to $\omega \mapsto T^*_\omega$, we obtain an operator $T^* \colon \Ind_\Omega F \to \Ind_\Omega E$ which is adjoint to $T$. 
\end{proof}
\end{prop}

\section{Crossed products by \'etale correspondences}\label{crossed product section}
Let $\Omega \colon G \to H$ be a Hausdorff \'etale correspondence and let $E \colon B \to C$ be an $H$-equivariant correspondence. In this section we will construct a correspondence 
\[\Omega \ltimes E \colon G \ltimes \Ind_\Omega B \to H \ltimes C\]
which we call the \textit{crossed product} of $E$ by $\Omega$. We build $\Omega \ltimes E$ from the space $\Gamma_c(\Omega, \sigma^* \E)$ of compactly supported sections of $\sigma^* \E \to \Omega$. Consider the following structure on $\Gamma_c(\Omega, \sigma^* \E)$:
\begin{itemize}
\item A $\Gamma_c(H, s^* \C)$-valued inner product. For each $\xi, \eta \in \Gamma_c(\Omega , \sigma^* \E)$, we define $\langle \xi,\eta \rangle \in \Gamma_c(H,s^* \C)$ for $h \in H$ by 
\begin{equation}\label{crossed product inner product formula}
 \langle \xi, \eta \rangle \colon  h \mapsto \sum_{ \omega \in \Omega_{r(h)} }  \langle h^{-1} \bcdot  \xi(\omega), \eta( \omega \bcdot h) \rangle \in C_{s(h)}.
\end{equation}
\item A right action $\Gamma_c(\Omega, \sigma^* \E) \rightacts \Gamma_c(H,s^* \C)$. For each $\xi \in \Gamma_c(\Omega , \sigma^* \E)$ and $c \in \Gamma_c(H,s^* \C)$ we define $\xi \bcdot c \in \Gamma_c(\Omega , \sigma^* \E)$ for $\omega \in \Omega$ by 
\begin{equation}\label{crossed product right action formula} 
\xi \bcdot c \colon  \omega \mapsto \sum_{h \in H^{\sigma(\omega)}} \left( h \bcdot  \xi  ( \omega \bcdot h )  \right) \bcdot  c (h^{-1} ) \in E_{\sigma(\omega)}. 
\end{equation}
\item A left action $\Gamma_c(G,s^* \Ind_\Omega \B) \acts \Gamma_c(\Omega, \sigma^* \E)$. For each $b \in \Gamma_c(G, s^* \Ind_\Omega \B)$ and $\xi \in \Gamma_c(\Omega , \sigma^* \E)$ we define $b \bcdot \xi \in \Gamma_c(\Omega,\sigma^* \E)$ for $\omega \in \Omega$ by
\begin{equation}\label{crossed product left action formula}
b \bcdot \xi\colon  \omega \mapsto \sum_{g \in G_{\rho(\omega)} }    (b(g^{-1})(g \bcdot \omega)) \bcdot \xi( g \bcdot \omega )  \in E_{\sigma(\omega)}. 
\end{equation}
\end{itemize}
The sections described in (\ref{crossed product inner product formula}), (\ref{crossed product right action formula}) and (\ref{crossed product left action formula}) are indeed compactly supported and continuous by applications of Lemma \ref{propersupportlemma} to the local homeomorphisms $(\omega,h) \mapsto h \colon \Omega \rtimes H \to H$, $(\omega,h) \mapsto \omega \bcdot h \colon \Omega \rtimes H \to \Omega$ and $(g,\omega) \mapsto g \bcdot \omega \colon G \ltimes \Omega \to \Omega$ respectively. Straightforward computations verify that the actions of $\Gamma_c(H, s^* \C)$ and $\Gamma_c(G, s^* \Ind_\Omega \B)$ on $\Gamma_c(\Omega, \sigma^* \E)$ are bilinear and respect the multiplication in these $*$-algebras, and that
\[\langle - , - \rangle \colon \Gamma_c(\Omega,\sigma^* \E) \times \Gamma_c(\Omega, \sigma^* \E) \to \Gamma_c(H, s^* \C)\]
is a Hermitian sesquilinear form on $\Gamma_c(\Omega,\sigma^* \E)$ compatible with its right action of $\Gamma_c(H, s^* \C)$. The construction of $\Omega \ltimes E \colon G \ltimes \Ind_\Omega B \to H \ltimes C$ breaks into two components: the Hilbert $H \ltimes C$-module $\Omega \ltimes E$ and the structure map $G \ltimes \Ind_\Omega B \to \cal L(\Omega \ltimes E)$. To complete $\Gamma_c(\Omega,\sigma^* \E)$ into a Hilbert $H \ltimes C$-module, we only need the data of the right action $\Omega \rightacts H$ and the $H$-Hilbert $C$-module $E$.

\begin{prop}[The crossed product Hilbert module]
Let $H$ be a Hausdorff \'etale groupoid, let $C$ be an $H$-$\Cast$-algebra, let $E$ be an $H$-Hilbert $C$-module and let $\Omega$ be a Hausdorff free proper \'etale right $H$-space with anchor map $\sigma \colon \Omega \to H^0$. Then there is a Hilbert $H \ltimes C$-module $\Omega \ltimes E$ containing a dense copy of $\Gamma_c(\Omega, \sigma^* \E)$ whose inner product and right module structure extend (\ref{crossed product inner product formula}) and (\ref{crossed product right action formula}).
\begin{proof}
We wish to complete $\Gamma_c(\Omega, \sigma^* \E)$ with respect to the norm 
\[ \lVert \xi \rVert \defeq \lVert \langle \xi, \xi \rangle \rVert_{H \ltimes C}^{\frac{1}{2}} \qquad \xi \in \Gamma_c(\Omega, \sigma^* \E) \] induced by $\langle - , - \rangle$ to obtain a Hilbert $H \ltimes C$-module $\Omega \ltimes E$. That this is in fact a norm hinges upon the positivity of $\langle - , - \rangle$. Let $\xi \in \Gamma_c(\Omega,\sigma^* \E)$. We aim to show that $\langle \xi , \xi \rangle$ is a positive element of $H \ltimes C$. There are finitely many slices $U_i \subseteq \Omega$ and $\xi_i \in \Gamma_c(U_i, \sigma^* \E)$ with $\lVert \xi_i \rVert_\infty \leq \lVert \xi \rVert_\infty$ such that $\xi = \sum_{i=1}^n \xi_i$. Consider the Hilbert $\Ind_\Omega C$-module $\Ind_\Omega E$. For each $i$, the inner product induces the supremum norm on $\Gamma_c(U_i, \sigma^* \E)$ so there is an isometric linear map $\Phi_i \colon \Gamma_c(U_i,\sigma^* \E) \to \Ind_\Omega E$ given by $H$-equivariant extension. Applying the left-sided version of \cite[Lemma 6.3]{MuhWil08a} to the left Hilbert $\mathcal K(\Ind_\Omega E)$-module $\Ind_\Omega E$, for any $\epsilon > 0$ there are finitely many elements $\nu_j \in \Ind_\Omega E$ such that 
\begin{equation*}
\sum_{j=1}^m \nu_j \bcdot \langle \nu_j , \Phi_k(\xi_k) \rangle \sim_\epsilon \Phi_k(\xi_k)
\end{equation*} 
for each $1 \leq k \leq n$. We use $x \sim_\epsilon y$ to mean that $d(x,y) < \epsilon$. For each $i$ the image of $\Phi_i$ is the set of equivariant sections supported in $U_i \bcdot H$ which have compact support when restricted to $U_i$, and $\Phi_i^{-1}$ is just the restriction map. In order to recover elements of $\Gamma_c(\Omega, \sigma^* \E)$ from the $\nu_j$ we take $\varphi_i \in C_c(q(U_i)) \subseteq C_c(\Omega/H)$ such that $\sum_{i=1}^n \lvert \varphi_i(z) \rvert^2 = 1$ for $z \in \bigcup_{k=1}^n q(\supp (\xi_k))$. We then set
\begin{equation*}
\eta_{i,j} := \Phi_i^{-1} (\varphi_i \bcdot \nu_j) \in \Gamma_c(U_i, \sigma^* \E),
\end{equation*} 
which maps $u \in U_i$ to $\varphi_i([u]_H) \nu_j(u)$. We claim that $\langle \xi, \xi \rangle$ is approximated by the positive operator $\sum_{i=1}^n \sum_{j=1}^m \langle \eta_{i,j} , \xi \rangle^* \langle \eta_{i,j} , \xi \rangle$. For each $i$, $j$ and $k$ the element $\eta_{i,j} \bcdot \langle \eta_{i,j} , \xi_k \rangle \in \Gamma_c(\Omega, \sigma^* \E)$ is supported on $U_k$ and $\Phi_k \left(\eta_{i,j} \bcdot \left\langle \eta_{i,j} , \xi_k \right\rangle \right) = (\varphi_i \bcdot \nu_j) \bcdot \left\langle \varphi_i \bcdot \nu_j , \Phi_k(\xi_k) \right\rangle$. Summing over $i$ and $j$, 
\begin{align*}
\sum_{i=1}^n \sum_{j=1}^m \Phi_k \left( \eta_{i,j} \bcdot \langle \eta_{i,j} , \xi_k \rangle \right) & = \sum_{i=1}^n \sum_{j=1}^m (\varphi_i \bcdot \nu_j) \bcdot \left\langle \varphi_i \bcdot \nu_j , \Phi_k(\xi_k) \right\rangle \\
& = \sum_{j=1}^m \nu_j \bcdot \langle \nu_j , \Phi_k(\xi_k) \rangle \\
& \sim_\epsilon \Phi_k(\xi_k),
\end{align*} 
and as $\Phi_k$ is isometric we get $\sum_{i=1}^n \sum_{j=1}^m \eta_{i,j} \bcdot \langle \eta_{i,j} , \xi_k \rangle \sim_\epsilon \xi_k$. For $e_i \in \Gamma_c(U_i, \sigma^* \E)$ and $e_j \in \Gamma_c(U_j, \sigma^* \E)$, the inner product $\langle e_i, e_j \rangle \in \Gamma_c(H, s^* \C)$ is supported on a slice and $\lVert \langle e_i, e_j \rangle \rVert \leq \lVert e_i \rVert \lVert e_j \rVert$. We may calculate
\begin{align*}
\langle \xi , \xi \rangle & = \sum_{k=1}^n \sum_{l=1}^n \langle \xi_k, \xi_l \rangle \\
& \sim_{n^2 \lVert \xi \rVert_\infty \epsilon} \sum_{k=1}^n \sum_{l=1}^n \left\langle \xi_k, \sum_{i=1}^n \sum_{j=1}^m \eta_{i,j} \bcdot \langle \eta_{i,j} , \xi_l \rangle \right\rangle \\
& =  \sum_{i=1}^n \sum_{j=1}^m \langle \eta_{i,j} , \xi \rangle^* \langle \eta_{i,j} , \xi \rangle.
\end{align*}
Sending $\epsilon$ to $0$, we conclude that $\langle \xi, \xi \rangle \geq 0$. 

If $\xi \in \Gamma_c(\Omega, \sigma^* \E)$ satisfies $\langle \xi, \xi \rangle = 0$, then for each $y \in H^0$ we have $\langle \xi , \xi \rangle (y) = \sum_{\omega \in \Omega_y} \langle \xi(\omega), \xi(\omega) \rangle = 0$, and therefore $\xi = 0$. This means that $\langle - , - \rangle$ induces a norm on $\Gamma_c(\Omega, \sigma^* \E)$. The completion $\Omega \ltimes E$ comes equipped with an inner product $\langle -,-\rangle \colon \Omega \ltimes E \times \Omega \ltimes E \to H \ltimes C$. By continuity, the right action $\Gamma_c(\Omega, \sigma^* \E) \times \Gamma_c(H, s^* \C) \to \Gamma_c(\Omega,\sigma^* \E)$ completes to a right action $\Omega \ltimes E \times H \ltimes C \to \Omega \ltimes E$. That $\langle - , - \rangle$ is compatible with the right action of $H \ltimes C$ follows from continuity. It follows that $\Omega \ltimes E$ is a Hilbert $H \ltimes C$-module. 
\end{proof}
\end{prop}

We note that for a Hausdorff \'etale correspondence $\Omega \colon G \to H$, the crossed product Hilbert $C^*(H)$-module $\Omega \ltimes C_0(H^0)$ is identical to the Hilbert $C^*(H)$-module $C^*(\Omega)$ constructed in \cite[Section 7]{AKM22}. However, we consider a left action by adjointable operators of the algebra $G \ltimes \Ind_\Omega C_0(H^0) \cong C^*(G \ltimes \Omega/H)$ rather than of $C^*(G)$.

\begin{prop}[The structure map of the crossed product correspondence]\label{crossed product correspondence}
Let $\Omega \colon G \to H$ be a Hausdorff \'etale correspondence and let $E \colon B \to C$ be an $H$-equivariant correspondence. Then the left action $\Gamma_c(G, s^* \Ind_\Omega \B) \acts \Gamma_c(\Omega, \sigma^* \E)$ described by formula (\ref{crossed product left action formula}) extends to a $*$-homomorphism $G \ltimes \Ind_\Omega B \to \cal L( \Omega \ltimes E)$. This gives $\Omega \ltimes E$ the structure of a correspondence from $G \ltimes \Ind_\Omega B$ to $H \ltimes C$, which is proper if $E$ is proper.
\begin{proof}
The left action $\Gamma_c(G, s^* \Ind_\Omega \B) \acts \Gamma_c(\Omega, \sigma^* \E)$ respects convolution and adjoints in $\Gamma_c(G, s^* \Ind_\Omega \B)$ in the sense that for each $a, b \in \Gamma_c(G, s^* \Ind_\Omega \B)$ and $\eta,\xi \in \Gamma_c(\Omega, \sigma^* \E)$ we have
\[  a \bcdot (b \bcdot \xi) = (a * b) \bcdot \xi \quad \text{and} \quad \langle a \bcdot \xi, \eta \rangle = \langle \xi, a^* \bcdot \eta \rangle. \]
To show that the image of the action map 
\[\Gamma_c(G, s^* \Ind_\Omega \B) \times \Gamma_c(\Omega, \sigma^* \E) \to \Gamma_c(\Omega, \sigma^* \E)\]
has dense span, it is enough to approximate $\xi \in \Gamma_c(U, \sigma^* \E)$ for an arbitrary slice $U \subseteq \Omega$, as these span $\Gamma_c(\Omega, \sigma^* \E)$. For each $\epsilon > 0$, we may take $\eta \in \Gamma_c(U, \sigma^* \B)$ such that $\lVert \eta \bcdot \xi - \xi \rVert_{\infty} < \epsilon$. By $H$-equivariant extension, we obtain an element $\tilde{\eta} \in \Ind_\Omega B$ which has compact support with respect to $\Omega/H$. The associated section $ \hat \eta \in \Gamma_0(G^0, \Ind_\Omega \B)$ therefore has compact support. By construction, $\lVert \hat \eta \bcdot \xi - \xi \rVert = \lVert \eta \bcdot \xi - \xi \rVert_{\infty} < \epsilon$.

Since the left action is non-degenerate and respects convolution and adjoints, it extends by Lemma \ref{technical lemma for crossed product} to a non-degenerate representation $\Gamma_c(G, s^* \Ind_\Omega \B) \to \cal L(\Omega \ltimes E)$. This in turn extends to a non-degenerate representation $G \ltimes \Ind_\Omega B \to \cal L(\Omega \ltimes E)$.

Now suppose that $E \colon B \to C$ is proper, so that $B$ acts by compact operators on $E$. The action of $\Ind_\Omega B$ on $\Omega \ltimes E$ therefore factors through the action of $\Ind_\Omega \cal K(E)$ on $\Omega \ltimes E$. For each slice $U \subseteq \Omega$, consider sections $\xi, \eta \in \Gamma_c(U, \sigma^* \E)$ and their $H$-equivariant extensions $\tilde \xi, \tilde \eta \in \Ind_\Omega E$. The operator $\Theta_{\tilde \xi, \tilde \eta} \in \cal K(\Ind_\Omega E) = \Ind_\Omega \cal K(E)$ acts on $\Omega \ltimes E$ as the compact operator $\Theta_{\xi, \eta}$. Across all slices, operators of the form $\Theta_{\tilde \xi, \tilde \eta} \in  \Ind_\Omega \cal K(E)$ have dense span, and we may conclude that $\Ind_\Omega B$ and therefore $G \ltimes \Ind_\Omega B$ acts on $\Omega \ltimes E$ by compact operators. 
\end{proof}
\end{prop}

\begin{example}\label{action correspondence crossed product}
Consider the action correspondence $\Omega \colon G \to G \ltimes Z$ associated to a Hausdorff étale groupoid $G$ and a locally compact Hausdorff $G$-space $Z$ with anchor map $\tau \colon Z \to G^0$. The induction functor $\Ind_\Omega$ coincides with the forgetful map $\tau_*$ which takes a $G \ltimes Z$-$\cs$-algebra $B$ and equips the underlying $\cs$-algebra with the natural action of $G$. In particular, the crossed products $G \ltimes \Ind_\Omega B$ and $(G \ltimes Z) \ltimes B$ may be identified. Through this identification, the crossed product $\Omega \ltimes E$ of a $G \ltimes Z$-equivariant correspondence $E \colon B \to C$ coincides with the usual crossed product $(G \ltimes Z) \ltimes E \colon (G \ltimes Z) \ltimes B \to (G \ltimes Z) \ltimes C$.
\end{example}

While the correspondence $\Omega \ltimes E \colon G \ltimes B \to H \ltimes C$ does not in general specialise to the correspondence $C^*(\Omega) \colon C^*(G) \to C^*(H)$ in the case without coefficients, the only difference is in the structure map, and we may recover $C^*(\Omega) \colon C^*(G) \to C^*(H)$ as follows. The $G$-equivariant map $\orho \colon \Omega/H \to G^0$ induces the structure map $C_0(G^0) \to M(C_0(\Omega/H))$ for a $G$-equivariant correspondence 
\[E_\Omega = C_0(\Omega/H) \colon C_0(G^0) \to \Ind_\Omega C_0(H^0),\] 
which is proper if $\Omega$ is proper.

\begin{prop}\label{recovering the correspondence}
Let $\Omega \colon G \to H$ be a Hausdorff \'etale correspondence. Then the composition 
\[ \Omega \ltimes C_0(H^0) \circ G \ltimes E_\Omega \colon G \ltimes C_0(G^0) \to H \ltimes C_0(H^0) \]
is isomorphic to the $\Cast$-correspondence $C^*(\Omega) \colon C^*(G) \to C^*(H)$.
\begin{proof}
The Hilbert $G \ltimes \Ind_\Omega C_0(H^0)$-module underlying the correspondence $G \ltimes E_\Omega$ is $G \ltimes \Ind_\Omega C_0(H^0)$ itself, after we identify $C_0(\Omega/H)$ with $\Ind_\Omega C_0(H^0)$. Thus composing $\Omega \ltimes C_0(H^0)$ with $G \ltimes E_\Omega$ does not change the underlying Hilbert $H \ltimes C_0(H^0)$-module, and the structure map $G \ltimes C_0(G^0) \to \mathcal L(\Omega \ltimes C_0(H^0))$ satisfies
\[ a \bcdot (b \bcdot e) = (a \bcdot b) \bcdot e \]
for each $a \in G \ltimes C_0(G^0)$, $b \in G \ltimes \Ind_\Omega C_0(H^0)$ and $e \in \Omega \ltimes C_0(H^0)$. Picking $a \in C_c(G)$ and $b \in C_0(\Omega/H)$ we may describe the element $a \bcdot b \in G \ltimes \Ind_\Omega C_0(H^0)$ as an element of $\Gamma_c(G, s^* \E_\Omega)$, which at $g \in G$ is the element $(a \bcdot b)(g) \in (\Ind_\Omega C_0(H^0))_{s(g)} \subseteq C_b(\Omega^{s(g)})$ given at $\omega \in \Omega^{s(g)}$ by \[(a \bcdot b)(g) \colon \omega \mapsto a(g) b([\omega]_H) .\]
We may therefore plug $a \in C_c(G)$, $b \in C_0(\Omega/H)$ and $e \in C_c(\Omega)$ into \eqref{crossed product left action formula} to compute $a \bcdot (b \bcdot e) \in C_c(\Omega)$ at $\omega \in \Omega$ as
\[ (a \bcdot (b \bcdot e)) (\omega) = ((a \bcdot b) \bcdot e) (\omega) = \sum_{g \in G_{\rho(\omega)}} a(g^{-1}) b([g \bcdot \omega ]_H) e(g \bcdot \omega). \]
This determines the structure map because elements of the form $b \bcdot e \in C_c(\Omega)$ are dense in $\Omega \ltimes C_0(H^0)$ and $C_c(G)$ is dense in $G \ltimes C_0(G^0)$. Since this agrees with the structure map of the correspondence $C^*(\Omega) \colon C^*(G) \to C^*(H)$ \cite[Equation (7.3)]{AKM22} under the identifications $C^*(G) = G \ltimes C_0(G^0)$, $C^*(\Omega) = \Omega \ltimes C_0(H^0)$ and $C^*(H) = H \ltimes C_0(H^0)$, the $\cs$-correspondences are isomorphic.
\end{proof}
\end{prop}

\begin{lem}\label{lemma for natural transformation}
Let $\Omega \colon G \to H$ be a Hausdorff \'etale correspondence and let $E \colon B \to C$ be an $H$-equivariant correspondence. Then there are isomorphisms
\begin{align*}
\Phi \colon & G \ltimes \Ind_\Omega E \otimes_{G \ltimes \Ind_\Omega C} \Omega \ltimes C  \cong \Omega \ltimes E \\
\Psi \colon & \Omega \ltimes B \otimes_{H \ltimes B} H \ltimes E  \cong \Omega \ltimes E
\end{align*}
of correspondences from $G \ltimes \Ind_\Omega B$ to $H \ltimes C$ such that for each 
\begin{align*}
\xi & \in \Gamma_c(G, s^* \Ind_\Omega \E) \subseteq G \ltimes \Ind_\Omega E, & \eta & \in \Gamma_c(\Omega, \sigma^* \C) \subseteq \Omega \ltimes C, \\
\mu & \in \Gamma_c(\Omega, \sigma^* \B) \subseteq \Omega \ltimes B, & \nu & \in \Gamma_c(H, s^* \E) \subseteq H \ltimes E, 
\end{align*}
we have $\Phi(\xi \otimes \eta) \in \Gamma_c(\Omega, \sigma^* \E) \subseteq \Omega \ltimes E$ and $\Psi(\mu \otimes \nu) \in \Gamma_c(\Omega, \sigma^* \E) \subseteq \Omega \ltimes E$ given by
\begin{align*}
\Phi(\xi \otimes \eta) \colon \omega & \mapsto \sum_{g \in G_{\rho(\omega)}} (\xi(g^{-1}) (g \bcdot \omega)) \bcdot \eta(g \bcdot \omega ), \\
\Psi(\mu \otimes \nu) \colon  \omega & \mapsto \sum_{h \in H^{\sigma(\omega)}} (h \bcdot \mu(\omega \bcdot h)) \bcdot \nu(h^{-1}).
\end{align*}
Furthermore, for each $T \in \cal L(\Ind_\Omega E)$ and $\lambda \in \Gamma_c(\Omega, \sigma^* \E) \subseteq \Omega \ltimes E$, the element $\Phi (G \ltimes T \otimes 1) \Phi^{-1}(\lambda) \in \Gamma_c(\Omega, \sigma^* \E) \subseteq \Omega \ltimes E$ is given by
\[ \Phi (G \ltimes T \otimes 1) \Phi^{-1}(\lambda) \colon \omega \mapsto T_\omega \lambda(\omega). \]
\begin{proof}
We define a bilinear map $u \colon \Gamma_c(G, s^* \Ind_\Omega \E) \times \Gamma_c(\Omega, \sigma^* \C) \to \Gamma_c(\Omega,\sigma^* \E)$ by setting
\[
u(\xi , \eta) \colon \omega \mapsto \sum_{g \in G_{\rho(\omega)}} (\xi(g^{-1}) (g \bcdot \omega)) \bcdot \eta(g \bcdot \omega ).
\]
for each $\xi \in \Gamma_c(G, s^* \Ind_\Omega \E)$ and $\eta \in \Gamma_c(\Omega, \sigma^* \C)$. The section $u(\xi, \eta)$ is continuous and compactly supported by an application of Lemma \ref{propersupportlemma} to the source map $(g,\omega) \mapsto \omega \colon G \ltimes \Omega \to \Omega$. The bilinear map $u$ preserves the inner product valued in $\Gamma_c(H,s^* \C) \subseteq H \ltimes C$, and intertwines the left action of $\Gamma_c(G, s^* \Ind_\Omega \B)$ and the right action of $\Gamma_c(H, s^* \C)$. It therefore extends to an inner product preserving map $\Phi \colon G \ltimes \Ind_\Omega E \otimes_{G \ltimes \Ind_\Omega C} \Omega \ltimes C \to \Omega \ltimes E$. If we restrict to a slice $U \subseteq \Omega$, the norm on $\Gamma_c(U,\sigma^* \E) \subseteq \Omega \ltimes E$ is the uniform norm. Given a slice $V \subseteq G$, $\xi \in \Gamma_c(V, s^* \Ind_\Omega \E)$ and $\eta \in \Gamma_c(U,\sigma^* \C)$, the section $\Phi(\xi \otimes \eta)$ is supported on the slice $V \bcdot U \subseteq \Omega$ and satisfies 
\[\Phi(\xi \otimes \eta)(g \bcdot \omega) = (\xi(g)(\omega)) \bcdot \eta(\omega) \]
for each $g \in V$ and $\omega \in U$. By Proposition \ref{density in Banach spaces}, sections of this form have dense span in $\Gamma_c(V \bcdot U, \sigma^* \E)$ because they are closed under the pointwise action of $C_0(V \bcdot U)$ and for each $g \bcdot \omega \in V \bcdot U$, the set 
\[\{(\xi(g)(\omega)) \bcdot \eta(\omega) \suchthat \xi \in \Gamma_c(V,s^* \Ind_\Omega \E), \; \eta \in \Gamma_c(U,\sigma^* \C) \}\] 
has dense span in $E_{\sigma(\omega)}$. We may conclude that $\Phi$ has dense image and is therefore an isomorphism of correspondences from $G \ltimes \Ind_\Omega B$ to $H \ltimes C$. The construction of $\Psi$ follows the exact same process.

For each $T \in \cal L(\Ind_\Omega E)$, the map sending $\lambda \in \Gamma_c(\Omega, \sigma^* \E)$ to the section $\omega \mapsto T_\omega \lambda(\omega) \in \Gamma_c(\Omega, \sigma^* \E)$ defines an action of $T$ on the dense subspace $\Gamma_c(\Omega, \sigma^* \E)$. This extends to a $*$-homomorphism $\psi \colon \cal L(\Ind_\Omega E) \to \cal L(\Omega \ltimes E)$ such that $\psi(T)(\lambda)(\omega) = T_\omega \lambda(\omega)$ for each $T \in \cal L(\Ind_\Omega E)$, $\lambda \in \Gamma_c(\Omega, \sigma^* \E)$ and $\omega \in \Omega$ through the use of Lemma \ref{technical lemma for crossed product} with the trivial groupoid. The operators $\psi(T) \circ \Phi$ and $\Phi \circ (G \ltimes T \otimes 1)$ coincide on simple tensors $\xi \otimes \eta$ with $\xi \in \Gamma_c(G, s^* \Ind_\Omega \E)$ and $\eta \in \Gamma_c(\Omega, \sigma^* \C)$, and therefore $\Phi(G \ltimes T \otimes 1 )\Phi^{-1} = \psi(T)$.
\end{proof}
\end{lem}

\section{The induction functor of an \'etale correspondence}\label{KK section}

To define the KK-theoretic induction functor $\Ind_\Omega \colon \KK^H \to \KK^G$ of a Hausdorff \'etale correspondence $\Omega \colon G \to H$, we make use of cutoff functions for the proper groupoid $\Omega \rtimes H$. This allows us to perform an averaging process to turn an operator into an $H$-equivariant operator.

\begin{defn}[Cutoff function]
Let $G$ be an \'etale groupoid. A \textit{cutoff function} for $G$ is a continuous function $c\colon  G^0 \to \bb R_{\geq 0}$ such that:
\begin{itemize}
\item for each $u \in G^0$, we have $\sum_{g \in G_u} c(r(g)) = 1 $,
\item the map $r \colon \supp (c \circ s) \to G^0$ is proper.
\end{itemize}
A \textit{cutoff function for an \'etale correspondence} $\Omega \colon G \to H$ is a cutoff function for $\Omega \rtimes H$. Because $s \colon G \to G^0$ is open we have $\supp (c \circ s) = s^{-1}(\supp c)$ and so the properness condition is equivalent to the compactness of $G^K_{\supp c}$ for each compact $K \subseteq G^0$.
\end{defn}

\begin{example}\label{étale homomorphism cutoff}
Let $\varphi \colon G \to H$ be an étale homomorphism of Hausdorff étale groupoids with associated étale correspondence $\Omega_\varphi = G^0 \times_{H^0} H$. Then $c \colon \Omega_\varphi \to \mathbb R_{\geq 0}$ given at $(x,h) \in \Omega_\varphi$ by 
\[ c(x,h) = \begin{cases} 1 & \text{if } h \in H^0, \\ 0 & \text{otherwise.}  \end{cases} \]
is a cutoff function.
\end{example}

\begin{prop}[Tu, Propositions 6.10 and 6.11 in \cite{Tu99a}]\label{cutoff existence}
If a Hausdorff \'etale groupoid $G$ admits a cutoff function, it is proper. Conversely, if $G$ is proper with a $\sigma$-compact orbit space, then it admits a cutoff function.
\end{prop}

\begin{defn}[Induced Fredholm operator]
Let $\Omega \colon G \to H$ be a second countable Hausdorff \'etale correspondence and let $(E,T) \in \bb E^H (B,C)$ be a Kasparov cycle. Let $c \colon  \Omega \to \bb R$ be a cutoff function. Recall that we may identify $\cal L(\Ind_\Omega E)$ with the $H$-equivariant sections in $\Gamma_b(\Omega, \cal L(\sigma^* \E))$. Under this identification, we define the \textit{induced Fredholm operator} $\Ind_{\Omega ,c} T \in \cal L ( \Ind_\Omega E)$ at $\omega \in \Omega$ by
\[ \Ind_{\Omega , c} T \colon  \omega \mapsto  \sum_{h \in H^{\sigma(\omega)}} c(\omega \bcdot h) (h \bcdot T_{s(h)}) \in \cal L(E_{\sigma(\omega)}). \]
This section $\Omega \to \cal L(\sigma^*  \E)$ is $H$-equivariant by construction, and bounded by the summation condition of $c$. To justify strict continuity, we may check that for each $\xi \in \Gamma_c(\Omega, \sigma^* \E)$, the sections 
\begin{align*}
\omega & \mapsto \Ind_{\Omega,c} T(\omega) \xi(\omega) & \omega & \mapsto \Ind_{\Omega,c} T(\omega)^* \xi(\omega)
\end{align*}
are continuous. Consider the continuous section
\begin{align*}
\Omega \rtimes H & \to r^* \sigma^* \E \\
(\omega,h) & \mapsto c(\omega \bcdot h)(h \bcdot T_{s(h)})(\xi(\omega)).
\end{align*}
Combining the properness condition of $c$ with the compact support of $\xi$, this has compact support. We may therefore apply Lemma \ref{propersupportlemma} to deduce that $\omega \mapsto \Ind_{\Omega,c} T(\omega) \xi(\omega) \colon \Omega \to \sigma^* \E$ is continuous. The same argument applies to the adjoint, so the induced Fredholm operator $\Ind_{\Omega,c} T \colon \Omega \to \cal L(\sigma^* \E)$ is strictly continuous and therefore defines an element of $\cal L(\Ind_\Omega E)$.
\end{defn}

\begin{lem}\label{inducedfredholmcompactperturbation}
Let $\Omega \colon G \to H$ be a second countable Hausdorff \'etale correspondence, let $c \colon \Omega \to \bb R$ be a cutoff function for $\Omega$ and let $(E,T) \in \bb E^H(B,C)$ be a Kasparov cycle with structure map $\varphi \colon B \to \cal L(E)$. Then for each $\omega \in \Omega$, the operator 
\[ \Ind_{\Omega,c} T(\omega) = \sum_{h \in H^{\sigma(\omega)}} c(\omega \bcdot h) (h \bcdot T_{s(h)}) \in \cal L(E_{\sigma(\omega)})\]
is a compact perturbation of $T_{\sigma(\omega)}$. Furthermore, the compact perturbation is continuous in $\omega$ in the sense that the maps
\begin{align*}
\sigma^* \B & \to \sigma^* \cal K(\E) & \sigma^* \B & \to \sigma^* \cal K(\E) \\
(\omega,b) & \mapsto \varphi_{\sigma(\omega)}(b)( \Ind_{\Omega,c} T(\omega) - T_{\sigma(\omega)} ) & (\omega,b) & \mapsto  ( \Ind_{\Omega,c} T(\omega) - T_{\sigma(\omega)} )  \varphi_{\sigma(\omega)}(b)
\end{align*}
are continuous.
\begin{proof}
That $\Ind_{\Omega,c} T(\omega)$ is a compact perturbation of $T_{\sigma(\omega)}$ for each $\omega \in \Omega$ follows straightforwardly from the Fredholm properties of $T$ and the summation condition of the cutoff function $c$. To check that the two maps $\sigma^* \B \to \sigma^* \cal K(\E)$ are continuous, we may by Proposition \ref{mapoutofbundle} check that any $b \in \Gamma_c(\Omega, \sigma^* \B)$ is mapped to a continuous section $\Omega \to \sigma^* \cal K(\E)$. The section
\begin{align*}
\Omega \rtimes H & \to r^* \sigma^* \cal K(\E) \\
(\omega,h) & \mapsto c(\omega \bcdot h) \varphi_{\sigma(\omega)}(b(\omega)) (h \bcdot T_{s(h)} - T_{\sigma(\omega)}) 
\end{align*}
is well-defined and continuous by almost invariance of $T$. It is compactly supported by the properness condition of $c$ combined with the compact support of $b \in \Gamma_c(\Omega, \sigma^* \B)$. Therefore by Lemma \ref{propersupportlemma}, the section
\begin{align*}
\Omega & \to \sigma^* \cal K(\E) \\
\omega & \mapsto \sum_{h \in H^{\sigma(\omega)}} c(\omega \bcdot h) \varphi_{\sigma(\omega)}(b(\omega)) (h \bcdot T_{s(h)} - T_{\sigma(\omega)}) \\
& = \varphi_{\sigma(\omega)}(b(\omega))( \Ind_{\Omega,c} T(\omega) - T_{\sigma(\omega)} )
\end{align*}
is well-defined, continuous and compactly supported. This shows that the first map $\sigma^* \B \to \sigma^* \cal K(\E)$ is continuous, and the same argument shows that the second is also continuous.
\end{proof}
\end{lem}
We check that the induced Fredholm operator really does define a Fredholm operator:
\begin{prop}[Induced Fredholm operator]
Let $\Omega \colon G \to H$ be a second countable Hausdorff \'etale correspondence, let $c \colon \Omega \to \bb R$ be a cutoff function for $\Omega$ and let $(E,T) \in \bb E^H(B,C)$ be a Kasparov cycle. Then $(\Ind_\Omega E, \Ind_{\Omega,c} T)$ is a $G$-equivariant Kasparov $\Ind_\Omega B$-$\Ind_\Omega C$ cycle.
\begin{proof}
Let $\varphi \colon  B \to \cal L (E)$ and $\hat{\varphi} \colon \Ind_\Omega B \to \cal L(\Ind_\Omega E)$ be the structure maps for $E \colon B \to C$ and $\Ind_\Omega E \colon \Ind_\Omega B \to \Ind_\Omega C$. We first check that for each $\xi \in \Ind_\Omega B$ the following operators in $\cal L(\Ind_\Omega E)$ are compact:
\begin{gather}
 \hat{\varphi}( \xi )  ((\Ind_{\Omega,c} T)^* - \Ind_{\Omega,c} T), \label{Fredholm one}\\
 \hat{\varphi}( \xi ) ((\Ind_{\Omega,c} T)^* \Ind_{\Omega,c} T - 1), \label{Fredholm two}\\
 [\Ind_{\Omega,c} T,\hat{\varphi}( \xi )]. \label{Fredholm three} 
\end{gather}
By Proposition \ref{compactoperatorsinducedmodule} we may identify $\cal K(\Ind_\Omega E)$ with $\Ind_\Omega \cal K(E)$, so it suffices to check that these operators define continuous sections $\Omega \to \sigma^* \cal K(\E)$. By the Fredholm properties of $T$, the following are continuous sections $\Omega \rtimes H \to r_{\Omega \rtimes H}^* \sigma^* \cal K(\E)$:
\begin{align*}
(\omega,h) & \mapsto c(\omega \bcdot h) (\hat \varphi (\xi)(\omega))(h \bcdot T^*_{s(h)} - h \bcdot T_{s(h)}) \\
(\omega,h) & \mapsto c(\omega \bcdot h)\left((\hat \varphi(\xi)(\omega))(h \bcdot T_{s(h)})  - (h \bcdot T_{s(h)})(\hat \varphi(\xi)(\omega))\right)
\end{align*}
These have proper support with respect to $(\omega,h) \mapsto \omega \colon \Omega \rtimes H \to \Omega$ by the properness condition of the cutoff function $c$. The evaluations of the operators (\ref{Fredholm one}) and (\ref{Fredholm three}) at $\omega \in \Omega$ may be written as the sum over $h \in H^{\sigma(\omega)}$ of the above operators at $(\omega,h)$. Therefore by Lemma \ref{propersupportlemma} these define continuous sections $\Omega \to \sigma^* \cal K(\E)$ and so (\ref{Fredholm one}) and (\ref{Fredholm three}) are compact. The operator (\ref{Fredholm two}) may be similarly shown to be compact by two applications of Lemma \ref{propersupportlemma}.

We now turn to almost invariance of $\Ind_{\Omega,c} T$. For each $\xi \in \Ind_\Omega B$, we aim to show that the section 
\begin{equation*}\label{almost invariance technical section one}   
g \mapsto \hat{\varphi}_{s(g)} ( \xi_{s(g)} ) \left( g^{-1} \bcdot (\Ind_{\Omega,c} T)_{r(g)} - (\Ind_{\Omega,c} T)_{s(g)} \right) \colon G \to s^*_G \cal L(\Ind_\Omega \E)
\end{equation*}
is a continuous section into the compact operators bundle $s^*_G \cal K(\Ind_\Omega \E)$. Through the identification of $\cal K(\Ind_\Omega \E)$ with $\Ind_\Omega \cal K(\E)$, it suffices by Lemma \ref{maps into induced bundle} to check that for each $\gamma \in C_c(G)$ the section $\eta \colon G \ltimes \Omega \to s^*_{G \ltimes \Omega} \cal L(\sigma^* \E)$ defined by
\begin{equation*}\label{almost invariance technical section}
\eta \colon (g, \omega) \mapsto \gamma(g) (\hat{\varphi}(\xi)(\omega)) ( \Ind_{\Omega,c} T(g \bcdot \omega) - \Ind_{\Omega,c} T(\omega) ) 
\end{equation*}
is a continuous section into the compact operators bundle $s^*_{G \ltimes \Omega} \sigma^* \cal K(\E)$ which vanishes at infinity with respect to $G \ltimes \Omega/H$. It vanishes at infinity with respect to $G \ltimes \Omega/H$ because $\xi \in \Ind_\Omega B$ vanishes at infinity with respect to $\Omega/H$ and $\gamma \in C_c(G)$ is compactly supported. Consider the local homeomorphism $\alpha \colon (g,\omega,h) \mapsto (g,\omega) \colon G \ltimes \Omega \rtimes H \to G \ltimes \Omega$ and the section $\nu \colon G \ltimes \Omega \rtimes H \to \alpha^* s^*_{G \ltimes \Omega} \sigma^* \cal K(\E)$ defined by
\begin{equation*}
\nu \colon (g,\omega,h) \mapsto \gamma(g) (c(g \bcdot \omega \bcdot h) - c(\omega \bcdot h))(\hat{\varphi}(\xi)(\omega)) (h \bcdot T_{s(h)} - T_{\sigma(\omega)}).
\end{equation*}
Almost invariance of $T$ ensures that $\nu$ is well-defined and continuous. Combining the compact support of $\gamma \in C_c(G)$ with the proper support of $c \circ s_{\Omega \rtimes H}$ with respect to $r_{\Omega \rtimes H} \colon \Omega \rtimes H \to \Omega$, the section $\nu$ has proper support with respect to $\alpha$. By Lemma \ref{propersupportlemma} applied to the local homeomorphism $\alpha \colon G \ltimes \Omega \rtimes H \to G \ltimes \Omega$ and the continuous section $\nu \colon G \ltimes \Omega \rtimes H \to \alpha^* s^*_{G \ltimes \Omega} \sigma^* \cal K(\E)$, the section $\nu_* = \eta \colon G \ltimes \Omega \to s^*_{G \ltimes \Omega} \sigma^* \cal K(\E)$ is continuous. This completes the proof of almost invariance, and the pair $(\Ind_\Omega E, \Ind_{\Omega,c} T)$ is therefore a Kasparov cycle. 
\end{proof}
\end{prop}
We can now define the induction functor of an \'etale correspondence.
\begin{defn}[The induction functor]\label{induction functor definition KK}
Let $\Omega \colon G \to H$ be a second countable Hausdorff \'etale correspondence. The \textit{induction functor} $\Ind_\Omega \colon \KK^H \to \KK^G$ is given by the following.
\begin{itemize}
\item An $H$-$\Cast$-algebra $B \in \KK^H$ is mapped to $\Ind_\Omega B \in \KK^G$.
\item A class $[E,T] \in \KK^H(B,C)$ is mapped to the class $[\Ind_\Omega E, \Ind_{\Omega,c} T] \in \KK^G(\Ind_\Omega B, \Ind_\Omega C)$, where $c \colon  \Omega \to \bb R$ is any cutoff function for $\Omega$.
\end{itemize}
Recall that by Proposition \ref{cutoff existence}, $\Omega$ has at least one cutoff function. Thankfully, it does not matter which one we pick.
\end{defn}

\begin{prop}[Well-definition of the induction functor]
Let $\Omega \colon G \to H$ be a second countable Hausdorff \'etale correspondence and let $B, C \in \KK^H$. For each cutoff function $c \colon \Omega \to \bb R$, the assignment $(E,T) \mapsto (\Ind_\Omega E, \Ind_{\Omega,c} T) \colon \bb E^H(B,C) \to \bb E^G(\Ind_\Omega B, \Ind_\Omega C)$ respects homotopy of Kasparov cycles. Furthermore, it is independent of the cutoff function.
\begin{proof}
Fix a cutoff function $c$ and suppose that $(E_0,T_0)$ and $(E_1,T_1)$ are homotopic Kasparov cycles via a Kasparov cycle $(F,S) \in \bb E^H(B, C([0,1],C))$. We aim to find a homotopy from $ (\Ind_\Omega E_0 , \Ind_{\Omega,c} T_0)$ to $ (\Ind_\Omega E_1 , \Ind_{\Omega,c} T_1)$. 
After identifying $\Ind_\Omega C([0,1],C)$ with $C([0,1],\Ind_\Omega C)$, we can take the homotopy to be the Kasparov cycle 
\[(\Ind_\Omega F, \Ind_{\Omega,c} S) \in \bb E^G(\Ind_\Omega B, \Ind_\Omega C([0,1],C)).\] 
Therefore $[E,T] \mapsto [\Ind_\Omega E, \Ind_{\Omega,c} T] \colon \KK^H(B,C) \to \KK^G(\Ind_\Omega B, \Ind_\Omega C)$ is well-defined for each cutoff function $c$.

Now suppose that $c_0$ and $c_1$ are two different cutoff functions. For any $\xi \in \Ind_\Omega B$, the operator $\xi(\omega) \bcdot (\Ind_{\Omega,{c_0}} T (\omega) - \Ind_{\Omega, c_1} T (\omega))$ is compact for each $\omega \in \Omega$ and varies continuously in $\omega$ by Lemma \ref{inducedfredholmcompactperturbation}. By Proposition \ref{compactoperatorsinducedmodule} this section defines a compact operator and so $\Ind_{\Omega , c_0}T$ and $\Ind_{\Omega, c_1}T$ are compact perturbations. The Kasparov cycles $(\Ind_\Omega E , \Ind_{\Omega, c_0} T)$ and $(\Ind_\Omega E , \Ind_{\Omega, c_1} T)$ are therefore homotopic, so $[\Ind_\Omega E, \Ind_{\Omega,c} T] \in \KK^G(\Ind_\Omega B, \Ind_\Omega C)$ is independent of the cutoff function $c$.
\end{proof}

\end{prop}

\begin{thm}[The induction functor]\label{KKinductionfunctor}
Let $\Omega \colon G \to H$ be a second countable Hausdorff \'etale correspondence. Then $\Ind_\Omega \colon \KK^H \to \KK^G$ defines homomorphisms of Kasparov groups and respects the Kasparov product and identity classes. In other words, $\Ind_\Omega$ is an additive functor.
\begin{proof}
The assignment $(E,T) \mapsto (\Ind_\Omega E, \Ind_{\Omega,c} T)$ preserves direct sums of Kasparov cycles, so defines a homomorphism of the Kasparov groups. The identity class at $B \in \KK^H$ is represented by the Kasparov cycle $(B,0) \in \bb E^H(B,B)$, which is mapped to the identity class $[\Ind_\Omega B,0] \in \KK^G(\Ind_\Omega B, \Ind_\Omega B)$. To show that $\Ind_\Omega$ respects the Kasparov product, we will show that if $(E_1,T_1) \in \bb E^H(A,B)$, $(E_2,T_2) \in \bb E^H(B,C)$ and $(E,T) \in (E_1,T_1) \hash_B (E_2,T_2)$, then 
\[(\Ind_\Omega E , \Ind_{\Omega, c} T) \in (\Ind_\Omega E_1 , \Ind_{\Omega, c} T_1) \hash_{\Ind_\Omega B} (\Ind_\Omega E_2, \Ind_{\Omega, c} T_2).\] 
Let $\varphi \colon  A \to \cal L(E)$ and $\varphi_1\colon  A \to \cal L(E_1)$ be the structure maps for $E$ and $E_1$. We need to show that:
\begin{itemize}
\item for each homogeneous $\xi_1 \in \Ind_\Omega E_1$, the operator $\theta_{\xi_1} \in \cal L( \Ind_\Omega E_2 , \Ind_\Omega E)$ given by $\xi_2 \mapsto \xi_1 \otimes \xi_2$ under the identification $\Ind_\Omega E_1 \otimes_{\Ind_\Omega B} \Ind_\Omega E_2 \cong \Ind_\Omega E$ satisfies 
\begin{nalign}\label{inducedconnectioncondition}
 \theta_{\xi_1} (\Ind_{\Omega,c} T_2) & - (-1)^{\deg (\xi_1)} (\Ind_{\Omega,c} T) \theta_{\xi_1}   & \in \cal K(\Ind_\Omega E_2, \Ind_\Omega E), \\
  \theta_{\xi_1} (\Ind_{\Omega,c} T_2)^* & - (-1)^{\deg (\xi_1)} (\Ind_{\Omega,c} T)^* \theta_{\xi_1}  & \in \cal K(\Ind_\Omega E_2, \Ind_\Omega E).
\end{nalign}
\item for each $\eta \in \Ind_\Omega A$, 
\begin{equation}\label{inducedpositivitycondition}
\Ind_\Omega(\varphi)(\eta) \left[ \Ind_{\Omega, c} T_1 \otimes 1 , \Ind_{\Omega ,c} T   \right] \Ind_\Omega(\varphi)(\eta^*) \geq 0 \; \text{mod} \; \cal K(\Ind_\Omega E).
\end{equation}
\end{itemize}
By Proposition \ref{compactoperatorsinducedmodule} we may identify $\cal K(\Ind_\Omega E_2, \Ind_\Omega E)$ with $\Ind_\Omega \cal K(E_2,E)$, so to prove that the operators in (\ref{inducedconnectioncondition}) are compact it suffices to show that they define continuous sections $\Omega \to \sigma^* \cal K(\E_2, \E)$. By the connection condition for $T$ and $T_2$, the section $\Omega \rtimes H \to r^* \sigma^* \cal K(\E)$ given by
\[(\omega,h) \mapsto c(\omega \bcdot h) \left( h \bcdot \left( \theta_{h^{-1} \bcdot \xi_1(\omega)} (T_2)_{s(h)} - (-1)^{\deg(\xi_1)} T_{s(h)} \theta_{h^{-1} \bcdot \xi_1(\omega)} \right)  \right) \] 
is well-defined and continuous. It has proper support with respect to $(\omega,h) \mapsto \omega \colon \Omega \rtimes H \to \Omega$ by the properness condition of the cutoff function $c$. By Lemma \ref{propersupportlemma} the section $\Omega \to \sigma^* \cal K(\E)$
\begin{align*}
\omega & \mapsto \sum_{h \in H^{\sigma(\omega)}} c(\omega \bcdot h) \left( h \bcdot \left( \theta_{h^{-1} \bcdot \xi_1(\omega)} (T_2)_{s(h)} - (-1)^{\deg(\xi_1)} T_{s(h)} \theta_{h^{-1} \bcdot \xi_1(\omega)} \right)  \right) \\
& = ( \theta_{\xi_1} (\Ind_{\Omega,c} T_2))(\omega) - (-1)^{\deg (\xi_1)} ((\Ind_{\Omega,c} T) \theta_{\xi_1})(\omega)
\end{align*}
is well-defined and continuous. The first operator in (\ref{inducedconnectioncondition}) is therefore compact and the same argument applied to $T^*$ rather than $T$ shows that the second is compact. To prove (\ref{inducedpositivitycondition}), consider the operator 
\[R := \Ind_\Omega(\varphi)(\eta) \left[ \Ind_{\Omega, c} T_1 \otimes 1 , \Ind_{\Omega ,c} T   \right] \Ind_\Omega(\varphi)(\eta^*).\]
To check that $R \gtrsim 0$, we claim that it is enough to find bounded continuous sections $k_U \colon U \to \sigma^* \cal K(\E)$ and $p_U \colon U \to \cal L(\sigma^* \E)$ for each $U$ in some open cover $\cal U$ of $\Omega$ such that for each $U \in \cal U$ and $\omega \in U$, we have $R(\omega) = p_U(\omega) + k_U(\omega)$ and $p_U(\omega) \geq 0$. This is because by Propositions \ref{compactoperatorsinducedmodule} and \ref{adjointableoperatorsinducedmodule}, the restrictions $\cal L(\Ind_\Omega E) \to \Gamma_b(U, \cal L(\sigma^* \E))$ induce an embedding of $\Cast$-algebras
\[ \overline{C_0(\Omega/H)\cal L(\Ind_\Omega E)}/\cal K(\Ind_\Omega E) \hookrightarrow \prod_{U \in \cal U} \Gamma_b(U, \cal L(\sigma^* \E))/\Gamma_b(U, \sigma^* \cal K(\E)). \]
Positivity of $R$ up to the compact operators may be checked in the larger $\Cast$-algebra, which is exactly our claimed sufficient condition. Let $\cal U$ be the open cover of slices $U \subseteq \Omega$. Let $a = \Ind_\Omega(\varphi)(\eta)$ so that $a_\omega = \varphi_{\sigma(\omega)}( \eta(\omega))$ for $\omega \in \Omega$ and let $U \in \cal U$. By the positivity condition for $T$ and $T_1$, the section $R_1 \defeq a [ \sigma^* T_1 \otimes 1, \sigma^* T] a^* \colon \Omega \to \cal L(\sigma^* \E)$
\begin{align*}
R_1 \colon \Omega & \to  \cal L(\sigma^* \E) \\
\omega & \mapsto a_\omega [(T_1)_{\sigma(\omega)} \otimes 1, T_{\sigma(\omega)}]a^*_\omega 
\end{align*}
may be written on $U$ as the sum $p_U + k$ of a bounded continuous section $k \colon U \to \sigma^* \cal K(\E)$ and a positive bounded continuous section $p_U \colon U \to \cal L(\sigma^* \E)$. Our aim now is to show that the section $R - R_1 \colon \Omega \to \cal L(\sigma^* \E)$ lands in the compacts bundle and is continuous into it. We will write $\sim$ to denote that the difference of two sections in $\Gamma_b(\Omega, \cal L(\sigma^* \E))$ lies in $\Gamma_b(\Omega, \sigma^* \cal K(\E))$. Let \[R_2 \defeq a [ \Ind_{\Omega,c} T_1 \otimes 1, \sigma^* T]a^* \colon \Omega \to \cal L(\sigma^* \E).\]
By Lemma \ref{inducedfredholmcompactperturbation}, the sections $\Ind_{\Omega,c} T$ and $\sigma^* T$ are compact perturbations, and therefore $R \sim R_2$. By the Fredholm properties of $T$, we have $[a^*, \sigma^* T] \sim 0$ and $[a, \sigma^* T] \sim 0$. Setting $b \defeq \Ind_{\Omega}(\varphi_1)(\eta)$, we may then calculate:
\begin{align*}
  R - R_1 & \sim R_2 - R_1 \\
   & =  a \left[  (\Ind_{\Omega,c} T_1 - \sigma^* T_1) \otimes 1 , \sigma^* T \right] a^* \\
  & \sim  \left[ a( (\Ind_{\Omega,c} T_1 - \sigma^* T_1) \otimes 1 ) a^* , \sigma^* T \right] \\
  & =  \left[ b(\Ind_{\Omega,c} T_1 - \sigma^* T_1)  b \otimes 1  , \sigma^* T \right]. 
\end{align*} 
The section $b(\Ind_{\Omega,c} T_1 - \sigma^* T_1) b \colon \Omega \to \cal L(\sigma^* \E_1)$ is pointwise compact and continuous into the compact operators bundle $\sigma^* \cal K(\E_1)$ again by Lemma \ref{inducedfredholmcompactperturbation}. For any $\nu \in \Gamma_b(\Omega, \sigma^* \cal K(\E_1))$ of degree $1$, the graded commutator $[\nu \otimes 1, \sigma^* T] = (\nu \otimes 1) \sigma^* T + \sigma^* T (\nu \otimes 1)$ satisfies $[\nu \otimes 1, \sigma^* T] \sim 0$. This holds because for any homogeneous $e,f \in \sigma^* E_1$ of opposite degree, it follows from $T$ being a $T_2$ connection for $E_1$ that
\begin{align*}
(\Theta_{e,f} \otimes 1) \sigma^* T & = \theta_e \theta_f^* (\sigma^* T) \\
& \sim (-1)^{\deg (f)} \theta_e (\sigma^* T_2) \theta_f^* \\
& \sim (-1)^{\deg (e) + \deg (f)} (\sigma^* T) \theta_e \theta_f^* \\
& =  - \sigma^* T (\Theta_{e,f} \otimes 1).
\end{align*}
As a result, $R - R_1 \sim 0$, and so $R$ may be written as the required sum $p_U +k_U$ on $U$, with $p_U \colon U \to \cal L(\sigma^* \E)$ positive and $k_U \colon U \to \sigma^* \cal K(\E)$ continuously compact. We may conclude that $R \gtrsim 0$ and therefore that $\Ind_{\Omega,c} T$ is a Kasparov product of $\Ind_{\Omega,c} T_1$ and $\Ind_{\Omega, c} T_2$.
\end{proof}
\end{thm}

\begin{rmk}\label{trick for functoriality remark}
To prove that $\Ind_\Omega \colon \KK^H \to \KK^G$ preserves the Kasparov product we could instead have appealed to Théor\`eme A.2.2 in the appendix of \cite{Lafforgue07} which states that any $\alpha \in \KK^H(B,C)$ can be decomposed into a finite string of elements induced by equivariant $*$-homomorphisms and inverses of invertible elements induced by equivariant $*$-homomorphisms. It would then suffice to check that $\Ind_\Omega(\alpha \circ \beta) = \Ind_\Omega(\alpha) \circ \Ind_\Omega(\beta)$ for elements $\alpha$ and $\beta$ in $\KK^H$ whenever either $\alpha$ or $\beta$ is induced by an equivariant $*$-homomorphism.
\end{rmk}

\begin{prop}[Naturality in KK of the crossed product by an étale correspondence]\label{induction natural transformation KK}
Let $\Omega \colon  G \to H$ be a second countable Hausdorff \'etale correspondence and let $[E,T] \in \KK^H(B,C)$ be a Kasparov cycle. Then the following diagram in $\KK$ commutes.
\begin{equation}\label{natural transformation diagram}
\begin{tikzcd}[column sep = large]
{G \ltimes \Ind_\Omega B} \arrow[rr, "{G \ltimes \Ind_\Omega [E,T]}"] \arrow[d, "{[\Omega \ltimes B , 0]}"] & & {G \ltimes \Ind_\Omega C } \arrow[d, "{[\Omega \ltimes C , 0]}"]  \\
{H \ltimes B} \arrow[rr, "{H \ltimes [E,T]}"] & & {H \ltimes C}
\end{tikzcd} 
\end{equation}
In other words, the assignment taking an $H$-$\Cast$-algebra $B$ to the Kasparov cycle $[\Omega \ltimes B, 0] \in \KK(G \ltimes \Ind_\Omega B, H \ltimes B)$ forms a natural transformation:
\[\begin{tikzcd}
	{\KK^H_{\phantom{H}}} &&& {\KK}
	\arrow[""{name=0, anchor=center, inner sep=0}, "{G \ltimes \Ind_\Omega -}", bend left = 30, from=1-1, to=1-4]
	\arrow[""{name=1, anchor=center, inner sep=0}, "{H \ltimes -}"', bend right = 30, from=1-1, to=1-4]
	\arrow["{\Omega \ltimes -}", shorten <=3pt, shorten >=3pt, Rightarrow, from=0, to=1]
\end{tikzcd}\]
\begin{proof}
Diagram (\ref{natural transformation diagram}) commutes at the level of correspondences through the isomorphisms 
\begin{align*}
\Phi \colon & G \ltimes \Ind_\Omega E \otimes_{G \ltimes \Ind_\Omega C} \Omega \ltimes C  \cong \Omega \ltimes E \\
\Psi \colon & \Omega \ltimes B \otimes_{H \ltimes B} H \ltimes E  \cong \Omega \ltimes E
\end{align*}
from Lemma \ref{lemma for natural transformation}. Furthermore, for each $R \in \cal L(\Ind_\Omega E)$, the operator 
\[\psi(R) \defeq \Phi( G \ltimes R \otimes 1) \Phi^{-1} \in \cal L(\Omega \ltimes E)\]
maps $\lambda \in \Gamma_c(\Omega, \sigma^* \E)$ to the section $\omega \mapsto R_\omega \lambda(\omega) \in \Gamma_c(\Omega, \sigma^* \E)$. Let $c \colon \Omega \to \bb R$ be a cutoff function for $\Omega \colon G \to H$. Through $\Phi$ it is immediate that 
\[ {[\Omega \ltimes E, \psi(\Ind_{\Omega,c} T)] = [G \ltimes \Ind_\Omega E , G \ltimes \Ind_{\Omega,c} T] \otimes_{G \ltimes \Ind_\Omega C} [\Omega \ltimes C , 0].} \] It remains to show that
\[ {[\Omega \ltimes E, \psi(\Ind_{\Omega,c} T)] = [\Omega \ltimes B ,0] \otimes_{H \ltimes B} [H \ltimes E , H \ltimes T].} \]
This boils down to checking that $\psi(\Ind_{\Omega,c} T)$ is an $H \ltimes T$-connection. For each $\mu \in \Omega \ltimes B$ let $T_\mu \in \cal L(H \ltimes E , \Omega \ltimes  E)$ be given by $\nu \mapsto \Psi(\mu \otimes \nu)$. Ultimately, we need to check that the operators
\begin{align*}
 S_\mu & := T_\mu (H \ltimes T) - \psi(\Ind_{\Omega,c} T) T_\mu & \in \cal L(H \ltimes E , \Omega \ltimes E) \\
 S'_\mu & := T_\mu (H \ltimes T)^* - \psi(\Ind_{\Omega,c} T)^* T_\mu & \in \cal L(H \ltimes E , \Omega \ltimes E)
\end{align*}
are compact. We will show that $S_\mu$ is compact by showing that $S^*_\mu S_\mu \in \cal L(H \ltimes E)$ is compact. To simplify notation we set $b_\omega \defeq \varphi_{\sigma(\omega)}(\mu(\omega))$ for $\omega \in \Omega$, where $\varphi \colon B \to \cal L(E)$ is the structure map for $E \colon B \to C$. Assume for now that $\mu \in \Gamma_c(\Omega, \sigma^* \B) \subseteq \Omega \ltimes B$, and consider the following:
\begin{itemize}
\item the element $\zeta \in \Gamma_c(\Omega, \sigma^* \cal K(\E))$ of the Hilbert $H \ltimes \cal K(E)$-module $\Omega \ltimes \cal K(E)$ given by
\[\zeta(\omega)  := b_\omega T_{\sigma(\omega)} - \Ind_{\Omega,c} T(\omega) b_\omega \in \cal K(E_{\sigma(\omega)}), \]
\item the structure map $\beta \colon H \ltimes \cal K(E) \to \cal K(H \ltimes E)$ of the crossed product of the proper $H$-equivariant correspondence $E \colon \mathcal K(E) \to C$.
\end{itemize}
Note that $\zeta$ is continuous by Lemma \ref{inducedfredholmcompactperturbation} and the Fredholm property of $T$ which implies that $\omega \mapsto [b_\omega, T_{\sigma(\omega)}] \colon \Omega \to \sigma^* \cal K(\E)$ is continuous. We claim that 
\[ \beta(\langle \zeta, \zeta \rangle) = S^*_\mu S_\mu. \]
While we may directly compute $S_\mu$ and $\beta(\langle \zeta, \zeta \rangle)$ on an element $\nu \in \Gamma_c(H, s^* \E) \subseteq H \ltimes E$, we do not have such a closed form expression for $S^*_\mu$. Instead this equality follows from a lengthy but straightforward computation that for each $\nu_1, \nu_2 \in \Gamma_c(H, s^* \E)$ we have 
\[ \langle S_\mu(\nu_1), S_\mu(\nu_2) \rangle = \langle \beta(\langle \zeta , \zeta \rangle ) (\nu_1), \nu_2 \rangle \in \Gamma_c(H, s^* B). \]
We conclude that $S_\mu^*S_\mu$ and therefore $S_\mu$ is compact when $\mu \in \Gamma_c(\Omega, \sigma^* \B)$. As $S_\mu$ varies continuously in $\mu$, it follows that $S_\mu$ is compact for all $\mu \in \Omega \ltimes B$. Similarly, $S_\mu'$ is compact, and so we have verified that $\psi(\Ind_{\Omega,c} T)$ is an $H \ltimes T$-connection. Thus $[\Omega \ltimes E, \psi(\Ind_{\Omega,c} T)]$ is the Kasparov product for both routes round diagram (\ref{natural transformation diagram}).
\end{proof}
\end{prop}

\begin{rmk}
As in Remark \ref{trick for functoriality remark} we could have appealed to \cite[Théor\`eme A.2.2]{Lafforgue07}, which implies that it suffices to check naturality of $\Omega \ltimes -$ on equivariant $*$-homomorphisms. This follows from Lemma \ref{lemma for natural transformation}.
\end{rmk}

\begin{defn}[Induction natural transformation]
Let $\Omega \colon G \to H$ be a second countable Hausdorff \'etale correspondence. The \textit{induction natural transformation} 
\[\alpha_\Omega \colon K_*(G \ltimes \Ind_\Omega - ) \Rightarrow K_*( H \ltimes -)\colon \KK^H \rightrightarrows \Ab_* \]
sends an $H$-$\Cast$-algebra $B$ to the map $K_*(\Omega \ltimes B) \colon K_*(G \ltimes \Ind_\Omega B) \to K_*(H \ltimes B)$.
\end{defn}

Suppose that $\Omega \colon G \to H$ is proper. Then $C^*(\Omega) \colon C^*(G) \to C^*(H)$ is proper and induces a map $K_*(C^*(\Omega)) \colon K_*(C^*(G)) \to K_*(C^*(H))$ in K-theory. Let $f_\Omega \in \KK^G(C_0(G^0), \Ind_\Omega C_0(H^0))$ be the element induced by the proper $G$-equivariant continuous map $\orho \colon \Omega/H \to G^0$. By Proposition \ref{recovering the correspondence}, we may recover $K_*(C^*(\Omega))$ in terms of $\Ind_\Omega$, $\alpha_\Omega$ and $f_\Omega$ with
\[K_*(C^*(\Omega)) = \alpha_\Omega(C_0(H^0)) \circ K_*(f_\Omega) \colon K_*(C^*(G)) \to K_*(C^*(H)).  \]

\section{Compatibility with composition of correspondences}\label{compatibility with composition section}

\begin{prop}\label{composition of correspondences and induction}
Let $\Omega \colon G \to H$ and $\Lambda \colon H \to K$ be Hausdorff \'etale correspondences and let $A$ be a $K$-Banach space. Then there is a $G$-equivariant isometric isomorphism $\varphi_A \colon \Ind_\Omega \Ind_\Lambda A \cong \Ind_{\Lambda \circ \Omega} A$ given by:
\begin{align*}
\varphi_A \colon  \Ind_\Omega \Ind_\Lambda A & \to \Ind_{\Lambda \circ \Omega} A \subseteq \Gamma_b(\Omega \times_H \Lambda , \sigma_{\Omega \times_H \Lambda}^* \A )\\
\xi & \mapsto \varphi_A(\xi) \\
& \quad \; \; [\omega , \lambda]_H \mapsto \xi(\omega)(\lambda) \\
\varphi_A^{-1} \colon  \Ind_{\Lambda \circ \Omega} A & \to \Ind_\Omega \Ind_\Lambda A \subseteq \Gamma_b(\Omega, \sigma_\Omega^* \Ind_\Lambda \A ) \\
\eta & \mapsto \varphi_A^{-1}(\eta) \\
& \quad \; \; \omega \mapsto \varphi_A^{-1}(\eta)(\omega) \\
& \qquad \quad \; \; \,  \lambda \mapsto \eta([\omega, \lambda]_H)
\end{align*}
\begin{proof}
It is clear that so long as both maps are well-defined, they are inverse to each other, isometric and $G$-equivariant. We justify well-definition by viewing both spaces inside a larger Banach space and showing that elements belonging to both spaces have dense span in each space. Let $f \colon \Omega \times_{H^0} \Lambda \to K^0$ be the map which takes $(\omega,\lambda) \in \Omega \times_{H^0} \Lambda$ to $\sigma_\Lambda(\lambda)$. We may view both $\Ind_\Omega \Ind_\Lambda A$ and $\Ind_{\Lambda \circ \Omega} A$ as subspaces of the Banach space $B$ of bounded sections $\Omega \times_{H^0} \Lambda \to f^* \A$. We define the linear isometries $\iota_1 \colon \Ind_\Omega \Ind_\Lambda A \hookrightarrow B$ and $\iota_2 \colon \Ind_{\Lambda \circ \Omega} A \hookrightarrow B$ for $\xi \in \Ind_\Omega \Ind_\Lambda A$, $\eta \in \Ind_{\Lambda \circ \Omega} A$ and $(\omega,\lambda) \in \Omega \times_{H^0} \Lambda$ by 
\begin{align*}
\iota_1(\xi)(\omega,\lambda) & \defeq \xi(\omega)(\lambda),  \\
\iota_2(\eta)(\omega,\lambda) & \defeq \eta([\omega,\lambda]_H).
\end{align*}

Consider slices $U \subseteq \Omega$ and $V \subseteq \Lambda$ with $\rho_\Lambda(V) \subseteq \sigma_\Omega(U)$. We construct the following diagram using the equivariant extensions from Proposition \ref{equivariant extensions}. 
\begin{equation}\label{diagram for proof of compatibility with composition}
\begin{tikzcd}[column sep = 20]
{\Gamma_c(V, \sigma_\Lambda^* \A)} \arrow[r, "\cong"', "\alpha"] \arrow[d, "\text{ext}", hook]  & {\Gamma_c(V \circ U, \sigma_{\Omega \times_H \Lambda}^* \A)} \arrow[r, "\text{ext}", hook]             & \Ind_{\Lambda \circ \Omega} A \arrow[rd, "\iota_2", hook] &   \\
{\Gamma_c(\sigma_\Omega(U), \Ind_\Lambda \A)} \arrow[r, "\cong"', "\beta"] & {\Gamma_c(U, \sigma_\Omega^* \Ind_\Lambda \A)} \arrow[r, "\text{ext}", hook] & \Ind_\Omega \Ind_\Lambda A \arrow[r, "\iota_1", hook]     & B
\end{tikzcd} 
\end{equation}
Note that for each $\xi \in \Gamma_c(V, \sigma_\Lambda^* \A)$, its $K$-equivariant extension $\text{ext}(\xi) \in \Ind_\Lambda A$ is supported on $\rho_\Lambda(\supp \xi) \subseteq \sigma_\Omega(U)$ as a section $H^0 \to \Ind_\Lambda \A$. The map $\alpha$ is defined through the homeomorphism $[u,v]_H \mapsto v \colon V \circ U \to V$, and for each $\xi \in \Gamma_c(V, \sigma_\Lambda^* \A)$ and $[\omega,\lambda]_H \in \Omega \times_H \Lambda$ we have \[\alpha(\xi)([\omega,\lambda]_H) = \begin{cases} \xi(h^{-1} \bcdot \lambda) & \text{there is} \; h \in H^{\sigma_\Omega(\omega)} \; \text{with} \; \omega \bcdot h \in U, \\ 0 & \text{otherwise}. \end{cases}\]
The map $\beta$ is defined through the homeomorphism $\sigma_\Omega \colon U \to \sigma_\Omega(U)$, and for $\xi \in \Gamma_c(V, \sigma_\Lambda^* \A)$ and $(\omega,\lambda) \in \Omega \times_{H^0} \Lambda$, we have 
\[ \beta(\text{ext}(\xi))(\omega)(\lambda) = \begin{cases} \sum_{k \in K^{\sigma_\Lambda(\lambda)}} k \bcdot \xi(\lambda \bcdot k) & \omega \in U, \\ 0 & \text{otherwise}. \end{cases} \]
We may compute both $\iota_2( \text{ext} ( \alpha(\xi) ))$ and $\iota_1(\text{ext} ( \beta(\text{ext} (\xi)) ))$ at $(\omega,\lambda) \in \Omega \times_{H^0} \Lambda$ to be
\[
\begin{cases} \sum_{k \in K^{\sigma_\Lambda(\lambda)}} k \bcdot \xi(h^{-1} \bcdot \lambda \bcdot k) & \text{there is} \; h \in H^{\sigma_\Omega(\omega)} \; \text{with} \; \omega \bcdot h \in U, \\ 0 & \text{otherwise}. \end{cases} 
\]
Diagram (\ref{diagram for proof of compatibility with composition}) therefore commutes. Slices of the form $V \circ U$ form an open cover of $\Omega \times_H \Lambda$, and so elements of the form $\text{ext}(\alpha(\xi))$ have dense span in $\Ind_{\Lambda \circ \Omega} A$ by Proposition \ref{equivariant extensions}. Over all slices $V \subseteq \Lambda$, $K$-equivariant extensions $\text{ext}(\xi)$ of $\xi \in \Gamma_c(V, \sigma_\Lambda^* \A)$ have dense span in $\Ind_\Lambda A$. Now suppose that $\sum_{i=1}^n \text{ext}(\xi_i)$ approximates $\eta \in \Gamma_c(\sigma_\Omega(U), \Ind_\Lambda \A) \subseteq \Ind_\Lambda A$. We may take $\gamma \in C_c(\sigma_\Omega(U))$ such that $\gamma = 1$ on $\supp \eta$ and define $\gamma \xi_i \in \Gamma_c(V \cap \rho_\Lambda^{-1}(\sigma_\Omega(U)), \sigma_\Lambda^* \A)$ by $v \mapsto \gamma(\rho_\Lambda(v)) \xi_i(v)$. Then $\sum_{i=1}^n \text{ext}(\gamma \xi_i)$ approximates $\eta$. This shows that for a fixed slice $U \subseteq \Omega$, the $H$-equivariant extensions $\text{ext}(\xi)$ of $\xi \in \Gamma_c(V, \sigma_\Lambda^* \A)$ over all slices $V \subseteq \Lambda$ with $\rho_\Lambda(V) \subseteq \sigma_\Omega(U)$ have dense span in $\Gamma_c(\sigma_\Omega(U), \Ind_\Lambda \A)$. Now letting $U$ vary, elements of the form $\text{ext}(\beta(\text{ext}(\xi)))$ have dense span in $\Ind_\Omega \Ind_\Lambda A$. We may conclude that $\iota_1(\Ind_\Omega \Ind_\Lambda A) = \iota_2(\Ind_{\Lambda \circ \Omega} A)$ and therefore $\varphi_A$ and $\varphi_A^{-1}$ define an isometric isomorphism between $\Ind_\Omega \Ind_\Lambda$ and $\Ind_{\Lambda \circ \Omega} A$.
\end{proof}
\end{prop}

\begin{prop}\label{product of cutoff functions}
Let $\Omega \colon G \to H$ and $\Lambda \colon H \to K$ be Hausdorff \'etale correspondences and let $c_\Omega \colon \Omega \to \bb R$ and $c_\Lambda \colon \Lambda \to \bb R$ be cutoff functions for $\Omega$ and $\Lambda$. Then there is a cutoff function $c \colon \Omega \times_H \Lambda \to \bb R$ for $\Lambda \circ \Omega$ given at $[\omega,\lambda]_H \in \Omega \times_H \Lambda$ by
\[ c([\omega,\lambda]_H) = \sum_{h \in H^{\sigma(\omega)}} c_\Omega(\omega \bcdot h) c_\Lambda(h^{-1} \bcdot \lambda).\]
\begin{proof}
We must check that $c$ is well-defined, continuous and satisfies the summation and properness conditions. The summation condition is clear and well-definition follows. To check continuity, we first check that the continuous map 
\[c_\Omega \times c_\Lambda \colon (\omega,\lambda) \mapsto c_\Omega(\omega) c_\Lambda(\lambda) \colon \Omega \times_{H^0} \Lambda \to \bb R \]
has proper support with respect to the local homeomorphism $q \colon  \Omega \times_{H^0} \Lambda \to \Omega \times_H \Lambda$. Let $V \subseteq \Omega \times_H \Lambda$ be compact, and let $W \subseteq \Omega \times_{H^0} \Lambda$ be a compact set with $q(W) = V$. We have $q^{-1}(V) = W \bcdot H$ and $ \supp ( c_\Omega \times c_\Lambda ) = \supp c_\Omega \times_{H^0} \supp c_\Lambda$, and therefore
\begin{align*}
 q^{-1}(V) \cap \supp (c_\Omega \times c_\Lambda) & = W \bcdot H \cap (\supp (c_\Omega) \times_{H^0} \supp (c_\Lambda)) \\
 & = \left( W \bcdot (\Omega \rtimes H)^{\pi_\Omega(W)}_{\supp (c_\Omega)} \right) \cap \Omega \times_{H^0} \supp (c_\Lambda). 
\end{align*}
The set $W \bcdot (\Omega \rtimes H)^{\pi_\Lambda(W)}_{\supp (c_\Omega)}$ is compact because $c_\Omega$ is a cutoff function and $W$ is compact. As $\Omega \times_{H^0} \supp (c_\Lambda)$ is closed, $q^{-1}(V) \cap \supp (c_\Omega \times c_\Lambda)$ is compact, and so $c_\Omega \times c_\Lambda$ has proper support with respect to $q \colon \Omega \times_{H^0} \Lambda \to \Omega \times_H \Lambda$. The function $c\colon  \Omega \times_H \Lambda \to \bb R$ is obtained by summing the values of $c_\Omega \times c_\Lambda$ over the fibres of the local homeomorphism $q \colon  \Omega \times_{H^0} \Lambda \to \Omega \times_H \Lambda$, and is therefore continuous by Lemma \ref{propersupportlemma}.

For the properness condition of $c \colon \Omega \times_H \Lambda \to \bb R$, we may show that $(\Omega \times_H \Lambda \rtimes K)^V_{\supp c}$ is compact, since $V \subseteq \Omega \times_H \Lambda$ was an arbitrary compact subset. Consider the set $W' \subseteq \Omega \times_{H^0} \Lambda$ given by
\[ W' = W \bcdot H \cap (\supp c_\Omega \times_{H^0} \Lambda) = W \bcdot (\Omega \rtimes H)^{\pi_\Omega(W)}_{\supp c_\Omega}, \]
which is compact because $W$ is compact and $c_\Omega$ is a cutoff function. Then there is a continuous surjection given by
\begin{align*}
W' \times_\Lambda (\Lambda \rtimes K)^{\pi_\Lambda(W')}_{\supp c_\Lambda} & \to (\Omega \times_H \Lambda \rtimes K)^V_{\supp c} \\
((\omega, \lambda), (\lambda,k)) & \mapsto ([\omega,\lambda]_H,k).
\end{align*}
As $c_\Lambda$ is a cutoff function it follows that $(\Omega \times_H \Lambda \rtimes K)^V_{\supp c}$ is compact, and therefore $c$ is a cutoff function.
\end{proof}
\end{prop}

\begin{prop}\label{composition of correspondences and KK induction}
Let $\Omega \colon G \to H$ and $\Lambda \colon H \to K$ be second countable Hausdorff \'etale correspondences. Then the $G$-equivariant $*$-isomorphisms $\varphi_A \colon \Ind_\Omega \Ind_\Lambda A \to \Ind_{\Lambda \circ \Omega} A$ from Proposition \ref{composition of correspondences and induction} associated to each $A \in \KK^K$ induce a natural isomorphism 
\[ \varphi_{\Omega, \Lambda} \colon \Ind_\Omega \Ind_\Lambda \cong \Ind_{\Lambda \circ \Omega} \colon \KK^K \rightrightarrows \KK^G. \]
\begin{proof}
Let $(E,T) \in \bb E^K(A,B)$ be a $K$-equivariant Kasparov cycle, let $c_\Omega \colon \Omega \to \bb R$ and $c_\Lambda \colon \Lambda \to \bb R$ be cutoff functions for $\Omega$ and $\Lambda$ and let $c \colon \Omega \times_H \Lambda \to \bb R$ be the cutoff function given at $[\omega,\lambda]_H \in \Omega \times_H \Lambda$ by
\[ c([\omega,\lambda]_H) = \sum_{h \in H^{\sigma(\omega)}} c_\Omega(\omega \bcdot h) c_\Lambda(h^{-1} \bcdot \lambda).\]
Then the induced elements of $\KK^G$ are given by
\begin{align*}
\Ind_\Omega \Ind_\Lambda [E,T] &= [\Ind_\Omega \Ind_\Lambda E, \Ind_{\Omega, c_{\Omega}} \Ind_{\Lambda, c_{\Lambda}} T], \\
\Ind_{\Lambda \circ \Omega} [E,T] &= [\Ind_{\Lambda \circ \Omega} E, \Ind_{\Lambda \circ \Omega, c} T].
\end{align*}
The identification $\varphi_E \colon \Ind_\Omega \Ind_\Lambda E \cong \Ind_{\Lambda \circ \Omega} E$ from Proposition \ref{composition of correspondences and induction} is compatible with the identifications $\varphi_A$ and $\varphi_B$. Furthermore, conjugation by $\varphi_E$ takes $\Ind_{\Lambda \circ \Omega, c} T$ to $\Ind_{\Omega,c_\Omega} \Ind_{\Lambda, c_\Lambda} T$. It follows that 
\[\KK^G(\varphi_A) \otimes_{\Ind_{\Lambda \circ \Omega} A} \Ind_{\Lambda \circ \Omega} [E,T] = \Ind_\Omega \Ind_\Lambda [E,T] \otimes_{\Ind_\Omega \Ind_\Lambda B} \KK^G(\varphi_B). \qedhere\]
\end{proof}
\end{prop}

\begin{rmk}
This recovers Le Gall's result that a Morita equivalence $G \sim_M H$ of Hausdorff étale groupoids induces an equivalence $\KK^G \simeq \KK^H$ of equivariant Kasparov categories \cite[Theor\`eme 7.5]{LeGall99}. Our equivalence is naturally isomorphic to Le Gall's equivalence. Moreover, for any étale Hilsum--Skandalis morphism $\Omega \colon G \to H$ (see Example \ref{HS morphism}), the induction functor $\Ind_\Omega \colon \KK^H \to \KK^G$ is naturally isomorphic to the functor Le Gall $\varphi^* \colon \KK^H \to \KK^G$ constructs from the associated generalised morphism $\varphi \colon G \to H$ \cite[Définition 7.1]{LeGall99}. By compatibility with composition (see also \cite[Theor\`eme 7.2]{LeGall99}) we need only consider this in the special case of an étale homomorphism $\varphi \colon G \to H$, because any étale Hilsum--Skandalis morphism may be decomposed into an étale homomorphism and the inverse of an étale homomorphism. The constructions $\Ind_{\Omega_\varphi}$ and $\varphi^*$ coincide through the use of the cutoff function for $\Omega_\varphi$ in Example \ref{étale homomorphism cutoff}.
\end{rmk}

\begin{prop}\label{composition of correspondences and crossed products}
Let $\Omega \colon G \to H$ and $\Lambda \colon H \to K$ be Hausdorff \'etale correspondences, let $A$ be a $K$-$\Cast$-algebra and consider the $G$-equivariant $*$-isomorphism $\varphi_A \colon \Ind_\Omega \Ind_\Lambda A \to \Ind_{\Lambda \circ \Omega} A$ from Proposition \ref{composition of correspondences and induction}. Then the following diagram of correspondences commutes up to isomorphism.
\[\begin{tikzcd}
G \ltimes \Ind_\Omega \Ind_\Lambda A \arrow[rr, "\Omega \ltimes \Ind_\Lambda A"] \arrow[d, "{(G \ltimes \varphi_A, G \ltimes \Ind_{\Lambda \circ \Omega} A)}"', "\cong"] &  & H \ltimes \Ind_\Lambda A \arrow[d, "\Lambda \ltimes A"] \\
G \ltimes \Ind_{\Lambda \circ \Omega} A \arrow[rr, "(\Lambda \circ \Omega) \ltimes A"]                  &  & K \ltimes A                                            
\end{tikzcd}\]
\begin{proof}
We first construct an isomorphism of Hilbert $K \ltimes A$-modules
\[\Phi \colon \Omega \ltimes \Ind_\Lambda A \otimes_{H \ltimes \Ind_\Lambda A} \Lambda \ltimes A  \cong (\Lambda \circ \Omega) \ltimes A\]
such that for each $\xi \in \Gamma_c(\Omega, \sigma_\Omega^* \Ind_\Lambda \A)$ and $\eta \in \Gamma_c(\Lambda, \sigma_\Lambda^* \A)$, we have $\Phi(\xi \otimes \eta) \in \Gamma_c(\Omega \times_H \Lambda, \sigma^* \A)$ given at $[\omega,\lambda]_H \in \Omega \times_H \Lambda$ by
\[\varphi(\xi, \eta) \colon [\omega,\lambda]_H \mapsto \sum_{h \in H^{\sigma(\omega)}} (\xi(\omega \bcdot h)(h^{-1} \bcdot \lambda)) \eta(h^{-1} \bcdot \lambda). \] 
The section $\varphi(\xi , \eta)$ is continuous and compactly supported by Lemma \ref{propersupportlemma} because it is obtained by summing the fibres over the quotient map $q \colon \Omega \times_{H^0} \Lambda \to \Omega \times_H \Lambda$ of the continuous compactly supported section 
\[(\omega,\lambda) \mapsto (\xi(\omega)(\lambda))\eta(\lambda) \colon \Omega \times_{H^0} \Lambda \to \pi_\Lambda^* \sigma_\Lambda^* \A.\]
The assignment $(\xi, \eta) \mapsto \varphi(\xi, \eta)$ respects the right action of $\Gamma_c(K, s^* \A)$ and respects the inner products in the sense that for each $\xi_1, \xi_2 \in \Gamma_c(\Omega, \sigma_\Omega^* \Ind_\Lambda \A)$ and $\eta_1, \eta_2 \in \Gamma_c(\Lambda, \sigma_\Lambda^* \A)$ we have 
\[ \langle \varphi(\xi_1 , \eta_1) , \varphi(\xi_2 , \eta_2) \rangle = \langle \xi_1 \otimes \eta_1 , \xi_2 \otimes \eta_2 \rangle. \]
We therefore obtain an inner product and right $K \ltimes A$-module preserving map $\Phi \colon \Omega \ltimes \Ind_\Lambda A \otimes_{H \ltimes \Ind_\Lambda A} \Lambda \ltimes A  \to (\Lambda \circ \Omega) \ltimes A$ such that $\Phi(\xi \otimes \eta) = \varphi(\xi, \eta)$ for each $\xi \in \Gamma_c(\Omega, \sigma_\Omega^* \Ind_\Lambda \A)$ and $\eta \in \Gamma_c(\Lambda, \sigma_\Lambda^* \A)$. To see that $\Phi$ is surjective, let $U \subseteq \Omega$ and $V \subseteq \Lambda$ be slices and let $\xi \in \Gamma_c(U, \sigma_\Omega^* \Ind_\Lambda \A)$ and $\eta \in \Gamma_c(V, \sigma_\Lambda^* \A)$. Then $\Phi(\xi \otimes \eta)$ is supported on $V \circ U \subseteq \Omega \times_H \Lambda$ and satisfies
\[ \Phi(\xi \otimes \eta) ([u,v]_H) = (\xi(u)(v))\eta(v) \]
for $u \in U$ and $v \in V$ with $\sigma_\Omega(u) = \rho_\Lambda(v)$. By Proposition \ref{density in Banach spaces}, sections of this form have dense span in $\Gamma_c(V \circ U, \sigma^* \A)$ because they are closed under the pointwise action of $C_0(V \circ U)$ and for each $[u,v]_H \in V \circ U$, the set 
\[ \{ (\xi(u)(v)) \eta(v) \suchthat \xi \in \Gamma_c(U, \sigma_\Omega^* \Ind_\Lambda \A), \; \eta \in \Gamma_c(V, \sigma_\Lambda^* \A) \} \]
has dense span in $A_{\sigma_\Lambda(v)}$. Letting $U$ and $V$ vary, the subspaces $\Gamma_c(V \circ U, \sigma^* \A)$ span $\Gamma_c(\Omega \times_H \Lambda, \sigma^* \A)$. Thus $\Phi$ has dense image and is therefore an isomorphism of Hilbert $K \ltimes A$-modules. The original diagram of correspondences commutes because a straightforward computation verifies that for $\xi \in \Gamma_c(\Omega, \sigma_\Omega^* \Ind_\Lambda \A)$, $\eta \in \Gamma_c(\Lambda, \sigma_\Lambda^* \A)$ and $\nu \in \Gamma_c(G, s^* \Ind_\Omega \Ind_\Lambda \A)$ we have
\[ (G \ltimes \varphi_A) (\nu) \bcdot \Phi(\xi \otimes \eta) = \Phi( \nu \bcdot \xi \otimes \eta). \qedhere \]
\end{proof}
\end{prop}

\begin{rmk}
Given an étale closed subgroupoid $H$ of an étale groupoid $G$, there is an étale correspondence $\Omega = G_{H^0} \colon G \to H$ which factors through the Morita equivalence $\Lambda = G_{H^0} \colon G \ltimes \Omega/H \to H$. For an $H$-$\cs$-algebra $B$, the induced algebra $\Ind^G_H B$ coincides with $\Ind_\Omega B$ (see \cite{Boenicke20}) and through Example \ref{action correspondence crossed product} we may make the identifications $G \ltimes \Ind_\Omega B \cong (G \ltimes \Omega/H) \ltimes \Ind_\Lambda B$ and $\Omega \ltimes B \cong \Lambda \ltimes B$. By Proposition \ref{composition of correspondences and crossed products}, the crossed product correspondence $\Lambda \ltimes B$ inherits invertibility from $\Lambda$ and so $\Omega \ltimes B$ becomes an imprimitivity bimodule from $G \ltimes \Ind^G_H B$ to $H \ltimes B$, recovering a version of Green's Imprimitivity Theorem \cite{Green78}.
\end{rmk}

\begin{cor}
Let $\Omega \colon G \to H$ and $\Lambda \colon H \to K$ be second countable Hausdorff \'etale correspondences and consider the natural isomorphism $\varphi_{\Omega, \Lambda} \colon \Ind_\Omega \Ind_\Lambda \cong \Ind_{\Lambda \circ \Omega} \colon \KK^K \rightrightarrows \KK^G$ from Proposition \ref{composition of correspondences and KK induction}. Then the following diagram of natural transformations commutes.
\[ \begin{tikzcd}
K_*(G \ltimes \Ind_\Omega \Ind_\Lambda -) \arrow[rr, "\alpha_\Omega \Ind_\Lambda", Rightarrow] \arrow[d, "{K_*(G \ltimes \varphi_{\Omega,\Lambda})}", "\cong"', Rightarrow] &  & K_*(H \ltimes \Ind_\Lambda -) \arrow[d, "\alpha_\Lambda", Rightarrow] \\
K_*(G \ltimes \Ind_{\Lambda \circ \Omega} -) \arrow[rr, "\alpha_{\Lambda \circ \Omega}", Rightarrow]                                                              &  & K_*(K \ltimes -)                                                     
\end{tikzcd} \]
\end{cor}

Given an étale correspondence $\Omega \colon G \to H$, we described how to recover the $\cs$-correspondence $C^*(\Omega) \colon C^*(G) \to C^*(H)$ in Proposition \ref{recovering the correspondence} from the crossed product $\Omega \ltimes -$ and the $G$-equivariant correspondence $E_\Omega = C_0(\Omega/H) \colon C_0(G^0) \to \Ind_\Omega C_0(H^0)$. When $\Omega$ is proper, $E_\Omega$ is induced by the equivariant proper continuous map $\orho_\Omega \colon \Omega/H \to G^0$. We show that the associated equivariant $*$-homomorphism $\orho_\Omega^* \colon C_0(G^0) \to \Ind_\Omega C_0(H^0)$ also respects the composition of étale correspondences.

\begin{prop}
Let $\Omega \colon G \to H$ and $\Lambda \colon H \to K$ be proper second countable Hausdorff \'etale correspondences and consider the $G$-equivariant $*$-isomorphism
\[\varphi_{C_0(K^0)} \colon \Ind_\Omega \Ind_\Lambda C_0(K^0) \to \Ind_{\Lambda \circ \Omega} C_0(K^0)\]
from Proposition \ref{composition of correspondences and induction}. Then the following diagram commutes.
\[ \begin{tikzcd}[column sep = large]
C_0(G^0) \arrow[rr, "\orho_\Omega^*"] \arrow[d, "\orho_{\Lambda \circ \Omega}^*"] &  & \Ind_\Omega C_0(H^0) \arrow[d, "\Ind_\Omega(\orho_\Lambda^*)"]                             \\
\Ind_{\Lambda \circ \Omega} C_0(K^0)                                  &  & \Ind_\Omega \Ind_\Lambda C_0(K^0) \arrow[ll, "{\varphi_{C_0(K^0)}}"', "\cong"]
\end{tikzcd} \]
\begin{proof}
Consider the identifications
\begin{align*}
\Ind_\Omega C_0(H^0) & \cong C_0(\Omega/H), \\
\Ind_\Omega \Ind_\Lambda C_0(K^0) & \cong C_0(\Omega \times_H (\Lambda /K)), \\
\Ind_{\Lambda \circ \Omega} C_0(K^0) & \cong C_0((\Omega \times_H \Lambda)/K).
\end{align*}
Under these identifications, $\Ind_\Omega(\orho^*_\Lambda)$ is induced by the proper map $[\omega,[\lambda]_K]_H \mapsto [\omega]_H \colon \Omega \times_H (\Lambda/K) \to \Omega/H$ and $\varphi_{C_0(K^0)}$ is induced by the homeomorphism \[ {[[\omega,\lambda]_H]_K \mapsto [\omega,[\lambda]_K]_H \colon (\Omega \times_H \Lambda)/K \to \Omega \times_H (\Lambda /K).}\] These $G$-equivariant proper maps along with $\orho_\Omega$ and $\orho_{\Lambda \circ \Omega}$ form a commutative diagram at the level of topological spaces, and so the diagram of $*$-homomorphisms commutes.
\end{proof}
\end{prop}

The results of this section imply that the decomposition
\[ K_*(C^*(\Omega)) = \alpha_\Omega(C_0(H^0)) \circ K_*(f_\Omega) \colon K_*(C^*(G)) \to K_*(C^*(H)) \]
described at the end of Section \ref{KK section} respects the composition of étale correspondences.

\end{document}